\documentclass[11pt]{article}

\usepackage{amsmath,amssymb,amsfonts,amscd,amsthm,graphics}
\usepackage[french]{babel}
\usepackage[matrix,arrow,curve]{xy}

\newcommand{\vol}{\text{\rm vol}} 
\newcommand{\Aire}{\text{\rm Aire}}
\newcommand{\Long}{\text{\rm Long}}
 
\newcommand{\DP}[2]{{\frac{\partial #1}{\partial #2}}} 
\newcommand{\tto}[1]{\ensuremath{\xrightarrow[#1]{}}} 
\newcommand{\til}[1]{\widetilde{#1}}
\newcommand{\tx}[1]{\text{\rm #1}}
   
\newcommand{\NN}{\mathbb{N}}  
   
\newcommand{\RR}{\mathbb{R}}   
\newcommand{\Ss}{\mathbb{S}}

\newcommand{\demo}{\noindent {\it \small D{\'e}monstration. }} 
 
\newtheorem{defi}{D{\'e}finition} 
\newtheorem{thm}{Th{\'e}or{\`e}me}    
\newtheorem{cor}{Corollaire}  
\newtheorem{prop}{Proposition} 
\newtheorem{lem}{Lemme}

\def\mc{\mathcal}

\title{Volume et courbure totale pour les hypersurfaces de l'espace
euclidien}  
  
\author{Alexandru Oancea\footnotemark[2]}
 
\date{\today} 

\begin{document}

\maketitle

\renewcommand{\thefootnote}{\fnsymbol{footnote}}
\addtocounter{footnote}{2}
\footnotetext{
CMAT, Ecole Polytechnique, 91128 Palaiseau Cedex, France 
et UMPA, Ecole Normale Sup{\'e}rieure de Lyon, 
46 All{\'e}e d'Italie, 69364 Lyon, France. Email: {\tt
oancea@math.polytechnique.fr}. Pendant l'{\'e}laboration de ce travail
l'auteur a aussi b{\'e}n{\'e}fici{\'e} du support mat{\'e}riel et logistique du 
  Laboratoire de
  Math{\'e}matiques de l'Universit{\'e} Paris-Sud, B{\^a}t. 425, 91405 Orsay Cedex, France.
}

\renewcommand{\thefootnote}{\arabic{footnote}}
\setcounter{footnote}{0}

\begin{abstract} 

Le th{\`e}me de l'article est fourni par la conjecture empirique
suivante: une hypersurface ferm{\'e}e 
$M^n \subset E^{n+1}$ de diam{\`e}tre ext{\'e}rieur
born{\'e} et de ``grand'' volume devrait n{\'e}cessairement {\^e}tre
``tr{\`e}s'' 
courb{\'e}e. On prend comme mesure pour la courbure le volume recouvert
sur la sph{\`e}re unit{\'e} $\Ss^n \subset E^{n+1}$ par l'application de
Gauss, {\`a} savoir la courbure totale $T(M)$. Notre {\'e}tude est
motiv{\'e}e par 
la preuve de Burago \cite{BZ} de l'in{\'e}galit{\'e} $\vol(M)\le R^2 T(M)$ 
pour 
$n=2$ et $M \subset B^{n+1}(R)$. 
Nous donnons un exemple explicite qui prouve qu'une
in{\'e}galit{\'e} de la forme $\vol(M) \le C R^n T(M)$ 
ne peut {\^e}tre vraie en toute g{\'e}n{\'e}ralit{\'e} 
en dimension $n \ge 3$. 
Nous
mettons toutefois en {\'e}vidence une condition suffisante sur la courbure
de Ricci sour laquelle l'in{\'e}galit{\'e} est v{\'e}rifi{\'e}e en dimension $3$.
Le
r{\'e}sultat principal de l'article est une in{\'e}galit{\'e} {\`a}
caract{\`e}re semi-local qui majore le volume d'un compact $K$ 
de $M$ par la courbure totale d'un ouvert $U$ qui le contient, sous
l'hypoth{\`e}se que la courbure de Gauss-Kronecker ne s'annule pas sur
$K$ et $U$. Pour $M \subset B^{n+1}(R)$ on prouve que 
$$\vol(K) \le \frac{C_n \, R^n}{d(K,\partial U)^n}T(U)$$
avec $d(K,\partial U)$ une distance mesur{\'e}e {\it sur la sph{\`e}re
$\Ss^n$} via 
l'application de Gauss. Le lieu des points de courbure nulle apparait
de cette fa{\c c}on comme une obstruction {\`a} avoir une in{\'e}galit{\'e}
globale comme celle de Burago en dimension sup{\'e}rieure. 

Nous obtenons diff{\'e}rentes variantes de l'in{\'e}galit{\'e} ayant un
caract{\`e}re isop{\'e}rim{\'e}trique et nous montrons qu'elles sont
optimales. Au passage, nous obtenons une in{\'e}galit{\'e}
isop{\'e}rim{\'e}trique ``inverse'' 
valable dans les espaces {\`a} courbure constante.
\end{abstract} 

\bigskip

\section*{Introduction} 

Une immersion isom{\'e}trique d'hypersurface 
$\varphi : M^n \longrightarrow E^{n+1}$ de classe $C^2$ est localement 
rigide en dehors du lieu des points de courbure nulle et on 
s'attend {\`a} ce que les propri{\'e}t{\'e}s m{\'e}triques de $M$ 
d{\'e}terminent de nombreuses propri{\'e}t{\'e}s de $\varphi$. Il y a
pourtant peu de r{\'e}sultats quantitatifs qui
estiment des caract{\'e}ristiques extrins{\'e}ques en termes d'objets
intrins{\'e}ques et qui soient valables en toute dimension (voir
\cite{BZ} pour plus de r{\'e}f{\'e}rences). 

Le r{\'e}sultat qui est {\`a} la source du pr{\'e}sent travail est une
in{\'e}galit{\'e} de Burago \cite{BZ} \S6.2 {\it en dimension $n=2$}. Si
$T(M)= \int_M |K| dV$ d{\'e}signe la courbure totale de
l'immersion, $R$ est le rayon d'une boule contenant $\varphi(M)$ et
${A}(M)$ d{\'e}signe l'aire de $M$ alors 
\begin{equation} \label{Burago depart}
{A}(M) \le R^2 T(M)
\end{equation}
avec {\'e}galit{\'e} uniquement pour le plongement isom{\'e}trique 
standard de la sph{\`e}re 
de rayon $R$. Nous nous sommes pos{\'e}s la question s'il y a un analogue de
(\ref{Burago depart}) en dimension sup{\'e}rieure. La r{\'e}ponse est
{\it n{\'e}gative} pour ce qui est d'une in{\'e}galit{\'e} globale du type
\begin{equation} \label{desir} 
\vol(M)\le C_nR^nT(M)
\end{equation} 
N{\'e}anmoins, nous obtenons des in{\'e}galit{\'e}s 
{\`a} caract{\`e}re isop{\'e}rim{\'e}trique {\it semi-local} faisant
intervenir la courbure totale d'un ouvert sur lequel l'application de
Gauss est non-d{\'e}g{\'e}n{\'e}r{\'e}e et exhibant le lieu des points de
courbure nulle comme une obstruction {\`a} l'existence d'in{\'e}galit{\'e}s
globales. 

L'article est structur{\'e} comme suit. Nous pr{\'e}sentons dans la
section \ref{le contre-exemple} un exemple explicite qui prouve la 
non-existence d'analogues de (\ref{Burago depart}) en dimension
sup{\'e}rieure. La section \ref{ineg globales} d{\'e}veloppe des formules
int{\'e}grales concernant les polyn{\^o}mes sym{\'e}triques des courbures
principales et la fonction support de l'immersion. Elles 
g{\'e}n{\'e}ralisent celles d{\'e}j{\`a} obtenues par
Minkowski et Kubota dans le cas convexe \cite{BF}, p.64 et par Hsiung
\cite{H} dans le cas non-convexe. Elles sont utilis{\'e}es d'une part pour
donner des conditions intrins{\'e}ques 
suffisantes pour l'in{\'e}galit{\'e} globale en dimension $n=3$ et
d'autre part pour raffiner par la suite l'in{\'e}galit{\'e}
semi-locale. Nous obtenons notamment la
\begin{prop} \label{Ricci} 
Soit $\varphi~: M^3 \longrightarrow \RR^4$ une immersion
  isom{\'e}trique de classe $C^2$ avec $M$ une vari{\'e}t{\'e} riemannienne
  lisse compacte sans bord. On suppose que $\varphi(M) \subset
  B^4(0,R)$. Si $M$ v{\'e}rifie $Ric \ge - \alpha / {R^2}$ avec 
  $0 < \alpha < 6$ alors 
  $$\vol(M) < \frac 6 {6-\alpha} R^3 T(M)$$
\end{prop} 
Il serait bien-sur int{\'e}ressant de trouver des contre-exemples {\`a}
l'in{\'e}galit{\'e} (\ref{desir}) satisfaisant $Ric \ge -{6}/{R^2} - \epsilon$.
Dans la section \ref{ineg locales} nous obtenons des in{\'e}galit{\'e}s
{\`a} caract{\`e}re semi-local. On 
utilise une m{\'e}thode d'estimation a-priori pour les
{\'e}quations de Monge-Amp{\`e}re r{\'e}elles d{\'e}velopp{\'e}e par 
Rauch et Taylor \cite{RT} et reprise par Aubin
\cite{Aubin}. L'op{\'e}rateur de type Monge-Amp{\`e}re qui apparait
naturellement dans notre contexte est $f \longmapsto
\det(H_f + f\tx{Id})$ agissant sur les fonctions d{\'e}finies sur
$\Ss^n$. Ceci sugg{\`e}re une relation - qui nous reste pour
l'instant cach{\'e}e - avec la th{\'e}orie spectrale du laplacien sur
$\Ss^n$. Notre r{\'e}sultat principal est le 
\begin{thm} On suppose que $M \subset B^{n+1}(0,R) \subset E^{n+1}$. 
Soit $U \subset M$  un ouvert sur lequel l'application de
Gauss est non-d{\'e}g{\'e}n{\'e}r{\'e}e. Il existe une constante $C_n$ qui
d{\'e}pend uniquement de la dimension telle que, si $K \subset U$ est un 
compact, on ait 
\begin{equation} \label{the semi best}
\vol(K) \le C_nR^n\frac1{\big(d_{\Ss^n}(K, \ \partial U) \big)^n} T(U)
\end{equation} 
et 
\begin{equation} \label{the best} 
\vol(K) \le  C_nR^n\frac1{\big(d_{\Ss^n}(K, \ \partial U) \big)^{n-1}}
T(U) + \frac{R}n\vol(\partial K) 
\end{equation}
o{\`u} $T(U)=\int_U|K|dV$ d{\'e}signe la courbure totale de $U$ et
$d_{\Ss^n}(K, \ \partial U)$ est la {\rm distance sph{\'e}rique locale}
entre $K$ et $\partial U$ (voir la d{\'e}finition \ref{la distance
locale}). 
\end{thm}

\noindent Le terme $d_{\Ss^n}(K, \ \partial U)$ tend vers z{\'e}ro
lorsque $U$ approche $K$. Cela fait que l'estimation (\ref{the best}) 
est d'autant meilleure que le compact $K$ est situ{\'e} {\`a} une grande
distance du 
lieu des points de courbure nulle. C'est dans ce sens qu'on
interpr{\`e}te le lieu des points de courbure nulle comme une
obstruction {\`a} l'in{\'e}galit{\'e} globale. L'obtention de 
(\ref{the best}) {\`a} partir de
(\ref{the semi best}) {\`a} l'aide des formules int{\'e}grales globales de 
type Minkowski occupe la section \ref{local to global}.
L'optimalit{\'e} de l'in{\'e}galit{\'e} (\ref{the best}) est
discut{\'e}e dans la section \ref{l'optimalite}
et l'article cl{\^o}t sur une in{\'e}galit{\'e} isop{\'e}rim{\'e}trique inverse
trait{\'e}e dans l'appendice.

\bigskip 

{\it Remerciements.} Ce travail constitue la premi{\`e}re partie de ma
th{\`e}se. Il a {\'e}t{\'e} pr{\'e}par{\'e} sous la direction bienveillante et
inspir{\'e}e de Claude Viterbo, qui saura trouver ici l'expression de ma
gratitude. J'ai profit{\'e} de remarques, suggestions et patience
d'{\'e}coute de la part de J{\'e}r{\^o}me Bertrand, Charles Boubel, 
Emmanuel Ferrand, 
Paul Gauduchon, Taoufik Hmidi,
Jean Lannes, Christophe Margerin, Nicolae Mihalache, 
Liviu Ornea, Ferit Ozt{\"u}rk, Pierre Pansu, Jean-Marc Schlenker 
et Costin V{\^a}lcu. Je leur suis reconnaissant.

\section{Un contre-exemple en dimension $n\ge 3$} \label{le contre-exemple}

Nous construisons pour tout $n \ge 3$ une suite $(M_k^n)_{_{k\ge 1}}$ 
d'hypersurfaces compactes plong\'ees dans 
$\RR^{n+1}$, diff{\'e}omorphes {\`a} $\Ss^n$ et de diam\`etre
ext\'erieur born\'e, telles que
$$\frac{\vol(M_k^n)}{T(M_k^n)}\stackrel{k}{\longrightarrow}\infty$$
o\`u $T(M) = \int_M |K| dV$ d\'esigne la courbure totale de $M$ et
$K=k_1\cdot k_2 \cdot \ldots \cdot k_n$ est la courbure de
Gauss-Kronecker, \'egale au produit des courbures principales de
$M$. La courbure totale mesure avec multiplicit\'es 
le volume recouvert sur $\Ss^n$ par
l'application de Gauss d\'efinie sur $M$. 

L'exemple prouve que l'in{\'e}galit{\'e} (\ref{desir}) ne peut {\^e}tre
v{\'e}rifi{\'e}e telle quelle en dimension  au moins $3$: tout ce que l'on
peut esp{\'e}rer de 
mieux est d'exhiber des classes de vari{\'e}t{\'e}s int{\'e}ressantes d'un
point de vue g{\'e}om{\'e}trique qui la satisfont. Ceci justifie en
particulier la proposition \ref{Ricci}. La
construction que nous pr{\'e}sentons nous a {\'e}t{\'e} sugg{\'e}r{\'e}e par
Jean-Marc Schlenker. 

Soit 
$$\chi~:\RR \longrightarrow \RR$$
une fonction $C^\infty$ telle que $\text{supp } \chi \subseteq
[-4\pi,4\pi]$, $0\le \chi \le 1$, $\chi \equiv 1 $ sur $[-2\pi,2\pi]$
et $\vert \chi '  \vert \le 1$. Nous d{\'e}finissons aussi 
$$\psi~: \RR^n     \longrightarrow \RR$$
$$\psi(x_1,\dots,x_n)=\chi(x_1)\dots\chi(x_n)$$
et on a $\vert \DP{\psi}{x_i} \vert \le 1$, 
$\text{supp } \psi \subseteq [-4\pi,4\pi]^n$, $\psi \equiv 1$ sur 
$[-2\pi,2\pi]^n$ et $ 0\le \psi \le 1$. 

Consid{\'e}rons la fonction
$$F_k~:\RR^n \longrightarrow \RR, \ k\in\NN^*, \ 0 < \alpha < 1$$
$$F_k(x)=\frac{\sin{kx_1}}{k^\alpha}\cdot \psi(x)$$

Nous d{\'e}finissons une hypersurface $M_k$ dans $\RR^{n+1}$ en tronquant le
graphe de $F_k$ au del{\`a} de $\vert x_i \vert =4\pi$ 
et en le refermant de fa{\c c}on lisse par une demi-sph{\`e}re. Pour tout $k$,
l'hypersurface $M_k$ est contenue dans la boule de rayon $8\pi$.

Prouvons que la famille $M_k$ fournit un contre-exemple {\`a} la
conjecture, c'est-{\`a}-dire 
$$\lim_{k\rightarrow\infty} \frac{\vol(M_k)}{T(M_k)}
\longrightarrow \infty$$

Tous les \'equivalents qui suivent seront consid\'er\'es pour $k
\rightarrow \infty$. Montrons d'abord que $\vol(M_k)\sim k^{1-\alpha}$. 
On a 
$$\DP{F_k}{x_1}=k^{1-\alpha} \cos{kx_1}\psi(x) + k^{-\alpha} \sin{k
  x_1} \DP{\psi}{x_1}$$ 
$$\DP{F_k}{x_i}=k^{-\alpha} \sin{k x_1} \DP{\psi}{x_i}, \ 2 \le i \le
n$$
et ceci entra{\^\i}ne 
$$\vol(\text{graph}({F_k}_{\vert_{[-4\pi,4\pi]^n}}))
= \int_{[-4\pi,4\pi]^n} \sqrt{1+\vert d F_k \vert^2} dV 
\sim k^{1-\alpha} \int_{-2\pi}^{2\pi} \vert \cos{kx_1} \vert dx_1 
\sim k^{1-\alpha}$$
d'o{\`u} $\vol(M_k)\sim k^{1-\alpha}$.

D'un autre c\^ot\'e, la normale unitaire au graphe de $F_k$ est 
$\nu = \frac{1}{\sqrt{1+\vert d F_k \vert^2} } (\DP{F_k}{x_1},\dots,
\DP{F_k}{x_n} , -1)$ et son
image est contenue dans un tube autour du grand cercle
d{\'e}termin{\'e} sur $\Ss^n$ par le plan $x_{n+1}Ox_1$. Le rayon de ce tube
est {\'e}quivalent {\`a} $\max_{ 2 \le i \le n} \vert \DP{F_k}{x_i} \vert
\sim k^{-\alpha} $. Puisque chaque point en dehors du grand cercle est
touch{\'e} au plus $k$ fois on d{\'e}duit 
$$T({M_k}) \le C k (k^{-\alpha})^{n-1} =
C k^{1-(n-1)\alpha}$$
o{\`u} $C$ est une constante qui ne d{\'e}pend que de $n$. 

Pour une autre constante $C'$ on aura 
$$\frac{\vol(M_k)}{T({M_k})} \ge C' k^{(n-2)\alpha} \
\xrightarrow{k \rightarrow \infty} \infty$$
et ceci montre la validit{\'e} de la construction.

Il est int{\'e}ressant de remarquer l'apparition de l'exposant
$n-2$. Ceci est \`a relier au fait que l'in\'egalit\'e (\ref{desir})
est valable en dimension $1 \le n \le 2$.

%%% Local Variables: 
%%% mode: latex
%%% TeX-master: "Gauss"
%%% End: 

%\setcounter{chapter}{1}
%\setcounter{section}{1}
%\setcounter{equation}{0}
\section{Formules int{\'e}grales} \label{ineg globales}

Nous pr{\'e}sentons des g{\'e}n{\'e}ralisations des formules int{\'e}grales 
de Minkowski \cite{BF,H} qui impliquent la fonction 
support d'une immersion et les polynomes sym{\'e}triques en les courbures
principales. On r{\'e}cup{\`e}re la preuve \cite{BZ} de (\ref{Burago
depart}), une preuve simple de (\ref{desir}) en dimension $n=1$ ou
lorsque $M$ est le bord d'un convexe (th{\'e}or{\`e}me d'Archim{\`e}de),
ainsi que des conditions suffisantes en dimension $n=3$.

\subsection{Notations}

Nous allons adopter un point de vue intrins{\'e}que: l'hypersurface  
sera repr{\'e}sent{\'e}e par une immmersion isom{\'e}trique $\varphi:M^n 
\longrightarrow \RR^{n+1}$, o{\`u}   
$M^n$ est une vari{\'e}t{\'e} riemannienne ferm{\'e}e orient{\'e}e de
dimension $n$.   
Soit $\nu:M\longrightarrow \RR^{n+1}$ le champ normal unitaire sur 
$(M,\varphi)$ et, suivant \cite{BZ}, mettons  
$$p=<\varphi,\nu>,\quad \quad q=\vert\varphi -<\varphi,\nu>\nu \vert$$ 
La fonction $p$ est 
appel{\'e}e ``fonction support'' de $(M,\varphi)$ et $q$ est  
la longueur de la composante de $\varphi$ tangente {\`a} $\varphi(M)$. 
Elles v{\'e}rifient $$p^2+q^2=\vert \varphi\vert ^2$$ 
 
Pour des vecteurs $u_1, \dots , 
u_{n+1} \in \RR^{n+1}$, on note $(u_1, \dots, u_{n+1})$ le volume du 
parall{\'e}lipip{\`e}de orient{\'e} qu'ils d{\'e}terminent.  
Pour des fonctions sur $M$ {\`a} valeurs vectorielles $u_1, \dots 
,u_k \in C^\infty(M,\RR^{n+1})$ et des  
$1$-formes sur $M$ {\`a} coefficients vectoriels  
$\alpha^{k+1}, \dots, \alpha^{n+1}\in \Omega^1(M,\RR^{n+1})$, o{\`u}  
$\alpha^j=\sum_1^n\alpha^j_idy^i$ et $\{ y^i\}$ est un syst{\`e}me local de 
coordonn{\'e}es sur $M$, on d{\'e}finit la $n-k+1$-forme sur $M$ 
$$(u_1,\dots , u_k,\alpha^{k+1}, \dots, \alpha^{n+1}) = 
\sum_{i_k,\dots, i_{n+1}}(u_1,\dots , u_k,\alpha^{k+1}_{i_{k+1}}, \dots, 
\alpha^{n+1}_{i_{n+1}}) \ dy^{i_{k+1}}\wedge\dots\wedge dy^{i_{n+1}}$$ 
 
Dans l'expression ci-dessus, l'interversion de 
deux vecteurs ou d'un vecteur et d'une $1$-forme change le 
signe. Par contre, l'interversion de deux $1$-formes ne change pas le 
signe.  
 
Si $k_1,\dots , k_n$ sont les courbures principales de $M$, on note 
$$S_k=\sum_{i_1<\dots<i_k} k_{i_1}\dots k_{i_k}, \quad 1\leq k \leq n$$ 
le $k$-i{\`e}me polyn{\^o}me sym{\'e}trique {\'e}l{\'e}mentaire en les
courbures principales. On pose $S_0=1$. On a
$S_k=\widetilde{S}_k / k!$ avec 
$$\widetilde{S}_k=\sum_{i_1,\dots,i_k \ \text{distincts}} 
k_{i_1}\dots k_{i_k}, \quad 1\leq k \leq n$$  
Posons aussi 
$$S_k^i=k_i \big( \sum_{i_1<\dots<i_{k-1} \ : \  i_j\neq  i}  
k_{i_1}\dots k_{i_{k-1}} \big) = 
k_i S_{k-1} (k_1,\dots,\hat{k}_i,\dots, k_n)$$ 
On remarque $S_2^i = \tx{Ric}(e_i,e_i)$, avec $e_i$ la direction de
courbure correspondant {\`a} $k_i$. 

\subsection{Formules de Minkowski
  g{\'e}n{\'e}ralis{\'e}es}\label{Minkowski general}
 
\begin{prop} L'identit{\'e} int{\'e}grale suivante est v{\'e}rifi{\'e}e sur $M$
\begin{equation}\label{generale} 
(n-k+1) \int_M p^{l-1}S_{k-1}dV \ = \  
k \int_M p^l S_k dV  -  
(l-1)\int_M p^{l-2} \big( \sum_{i=1}^n 
S_k^i <\varphi,e_i>^2 \big)dV  
\end{equation}  
\hfill{$1\le k \le n, \ l \ge 1$} 

o{\`u} les $e_i$ d{\'e}signent les directions de courbure sur $M$. 
\end{prop} 

\demo. Nous allons calculer pour $1\leq k\leq n$ et $l\ge 1$ la 
diff{\'e}rentielle ext{\'e}rieure   
\begin{eqnarray*} 
d\big(p^{l-1}(\varphi,\nu,\underbrace{d\nu,\dots,d\nu}_{k-1}, 
\underbrace{d\varphi,\dots,d\varphi}_{n-k })\big) & = &  
(l-1)p^{l-2} dp \wedge (\varphi,\nu,\underbrace{d\nu,\dots,d\nu}_{k-1}, 
\underbrace{d\varphi,\dots,d\varphi}_{n-k }) \ - \\ 
& & - \ p^{l-1}(\nu,\underbrace{d\nu,\dots,d\nu}_{k-1}, 
\underbrace{d\varphi,\dots,d\varphi}_{n-k+1}) \ + \\ 
& & + \ p^{l-1}(\varphi,\underbrace{d\nu,\dots,d\nu}_k, 
\underbrace{d\varphi,\dots,d\varphi}_{n-k }) 
\end{eqnarray*} 
Autour de tout point $x\in M$ on peut choisir un syst{\`e}me de  
coordonn{\'e}es $\{y^i\}$ tel que, {\it au point $x$}, les 
$\{\DP{}{y^i}\}_{i=1\dots n}$ forment une base orthonorm{\'e}e et 
orient{\'e}e avec
$e_i=\varphi_*\DP{}{y^i}$ des directions de courbure 
au point $\varphi(x)$.  Les suivantes {\'e}galit{\'e}s sont valables au
point $x$~:  
\begin{eqnarray}\label{coords} 
d\nu=\sum k_ie_idy^i & \\ 
d\varphi=\sum e_idy^i & \nonumber \\  
<e_i,e_j>=\delta_{ij}, & \quad e_1 \times \dots \times e_n = \nu \nonumber  
\end{eqnarray} 
Dans un tel syst{\`e}me de coordonn{\'e}es on a
\begin{eqnarray}\label{nu} 
(\nu,\underbrace{d\nu,\dots,d\nu}_{k-1},\underbrace{d\varphi,\dots,  
d\varphi}_{n-k+1 }) & = & \sum_{\sigma \in S_n}(\nu,e_{\sigma(1)},\dots, 
e_{\sigma(n)})k_{\sigma(1)}\dots k_{\sigma(k-1)}dy^{\sigma(1)}\dots 
dy^{\sigma(n)} \nonumber \\ 
& = & (\sum_{\sigma \in S_n}k_{\sigma(1)}\dots 
k_{\sigma(k-1)})(\nu,e_1,\dots,e_n)dy^1\dots dy^n \nonumber \\ 
& = & (-1)^n(n-k+1)!(k-1)!S_{k-1}dV  
\end{eqnarray} 
\begin{eqnarray}\label{phi} 
(\varphi,\underbrace{d\nu,\dots,d\nu}_k,\underbrace{d\varphi,\dots,  
d\varphi}_{n-k }) & = & \sum_{\sigma \in S_n} (\varphi,e_{\sigma(1)},\dots, 
e_{\sigma(n)})k_{\sigma(1)}\dots k_{\sigma(k)}dy^{\sigma(1)}\dots 
dy^{\sigma(n)} \nonumber \\ 
& = &  (n-k)! k! S_k (\varphi,e_1,\dots , e_n) dV \nonumber \\ 
& = & (-1)^n (n-k)! k! S_k p dV   
\end{eqnarray} 
On a utilis{\'e} les identit{\'e}s  
$$(\nu,e_1,\dots,e_n)=(-1)^n(e_1,\dots,e_n,\nu)=(-1)^n<e_1\times\dots\times  
e_n,\nu>=(-1)^n$$ 
$$(\varphi,e_1,\dots,e_n)=(-1)^n(e_1,\dots,e_n,\varphi)=(-1)^n<e_1\times\dots\times  
e_n,\varphi>=(-1)^n p$$ 
Pour ce qui est du premier terme on a 
$$dp = d<\varphi,\nu> = \sum k_i<\varphi,e_i>dy^i$$ 
\begin{eqnarray*} 
(\varphi,\nu,\underbrace{d\nu,\dots,d\nu}_{k-1}, 
\underbrace{d\varphi,\dots,d\varphi}_{n-k }) = \hspace{4cm} \\  
\sum_{i=1}^n \sum_{\sigma \in S_n \ : \ \sigma(n)=i} 
(\varphi,\nu,e_{\sigma(1)},\dots,e_{\sigma(n-1)}) 
k_{\sigma(1)}\dots k_{\sigma(k-1)}dy^{\sigma(1)}\dots 
dy^{\sigma(n-1)} = \\ 
\sum_{i=1}^n 
\big(\sum_{\sigma \in S_n \ : \ \sigma(n)=i}k_{\sigma(1)}\dots 
k_{\sigma(k-1)}\big) (\varphi,\nu,e_1,\dots,\hat{e}_i,\dots,e_n)dy^1\dots  
d\hat{y}^i\dots dy^n  = \\ 
(n-k)! (k-1)! \sum_{i=1}^n  
\big( \sum_{i_1<\dots<i_{k-1} \ : \ i_j \neq i}  
k_{i_1}\dots k_{i_{k-1}} \big)  
(-1)^{n+i} <\varphi,e_i> dy^1\dots d\hat{y}^i\dots dy^n  
\end{eqnarray*} 
On obtient  
\begin{eqnarray}\label{phinu} 
dp\wedge (\varphi,\nu,\underbrace{d\nu,\dots,d\nu}_{k-1}, 
\underbrace{d\varphi,\dots,d\varphi}_{n-k }) = \hspace{4cm} \\ 
(-1)^{n-1} (n-k)! (k-1)! 
\Big( \sum_{i=1}^n k_i \big( \sum_{i_1<\dots<i_{k-1}\ : \ i_j \neq  i}  
k_{i_1}\dots k_{i_{k-1}} \big) <\varphi,e_i>^2 \Big)dV \ 
= \nonumber \\ 
(-1)^{n-1} (n-k)! (k-1)! \big( \sum_{i=1}^n S_k^i<\varphi,e_i>^2 
\big)dV \hspace{3,5cm} \nonumber  
\end{eqnarray} 
et, finalement,  
\begin{eqnarray*} 
(l-1) (n-k)! (k-1)! p^{l-2} \big( \sum_{i=1}^n 
S_k^i<\varphi,e_i>^2 \big)dV + 
(n-k+1)!(k-1)!p^{l-1}S_{k-1}dV - \\ 
- (n-k)! k! p^l S_k dV  = \ \text{forme exacte} \hspace{5cm} 
\end{eqnarray*} 
Il suffit maintenant d'int{\'e}grer cette derni{\`e}re identit{\'e} sur
$M$. 

\hfill{$\square$}

\subsection{Applications}\label{Applications}

a) Formules de Minkowski. Prendre $l=1$ pour obtenir
\begin{equation}\label{minkowski} 
(n-k+1)\int_M S_{k-1} dV = k \int_M p S_k dV,\quad 1\leq k \leq n 
\end{equation} 

b) L'in{\'e}galit{\'e} (\ref{desir}) pour $n=1$. Prendre $k=n=1$ dans les
formules de Minkowski pour obtenir, lorsque $C$ est une courbe
immerg{\'e}e dans le plan 
$$l(C) = \int_C pK dl$$
et ceci implique tout de suite $\l(C) \le R T(C)$, avec {\'e}galit{\'e} si 
et seulement si $C$ est le plongement standard du cercle de rayon
$R$. 

c) Expression du volume de $M$ en termes de courbure moyenne. Prendre
$k=1$ dans (\ref{minkowski}) et obtenir 
\begin{equation}\label{courbure moyenne}
\vol(M) = \frac1n \int_M pH dV
\end{equation} 

d) Expression du volume de $M$ en termes de courbure de Gauss. 
Prendre $l=k$ dans (\ref{generale}) pour obtenir 
\begin{equation}\label{l=k}
(n-k+1)\int_M p^{k-1} S_{k-1} dV = k\int_M p^k S_k dV - 
(k-1)\int_M p^{k-2} (\sum_{i=1}^n S_k^i <\varphi, e_i> ^2) dV
\end{equation} 
En additionnant les identit{\'e}s (\ref{l=k}) pour $1\le k \le n$ nous
trouvons 
\begin{equation} \label{formule globale generale}
\text{vol}(M)= \int_M \big( p^n - \frac{n-1}{n} p^{n-2} q^2 \big) K
dV - \sum_{k=2}^{n-1} C_{n,k} \int_M p^{k-2} (\sum_{i=1}^n S_k^i
<\varphi, e_i>^2 ) dV
\end{equation} 
avec $\displaystyle C_{n,k}= \frac{(k-1)!}{n(n-1)\dots(n-k+2)} \cdot
\frac{k-1}{n-k+1}$. 

e) Preuve de (\ref{desir}) lorsque $M$ est le bord d'un
convexe (th{\'e}or{\`e}me d'Archim{\`e}de). Toutes les courbures principales
sont alors positives et on 
obtient $\vol(M) \le \int_M (p^n - \frac{n-1}{n}p^{n-2}q^2)K dV \le
R^n T(M) = R^n \vol(\Ss^n)$, avec {\'e}galit{\'e} si et seulement si $M$
est la sph{\`e}re de rayon $R$. 

f) Preuve de (\ref{Burago depart}) selon \cite{BZ}. Prendre $n=2$
dans (\ref{formule globale generale}) pour obtenir 
$$A(M)=\int_M (p^2 - \frac12 q^2) K dS$$ 
Lorsque $\varphi(M)\subset B(0,R)$ on a $p^2 +q^2 \le R^2$ et cela
implique tout de suite $A(M) \le R^2T(M)$, avec {\'e}galit{\'e} si et
seulement si $\varphi$ est le plongement standard de la sph{\`e}re de
rayon $R$ dans $E^3$. 

g) Condition suffisante pour $n=3$. D'apr{\`e}s (\ref{formule globale
generale}) on a  
\begin{equation}\label{n=3}
\text{vol}(M)=\int_M (p^3 - \frac23 pq^2) K dV - \frac16 \int_M
\sum_{i=1}^3 Ric(e_i,e_i) <\varphi,e_i>^2 dV
\end{equation}
L'in{\'e}galit{\'e} (\ref{desir}) sera vraie sous la forme $\vol(M) \le R^3
 T(M)$ si $M$ v{\'e}rifie $\tx{Ric} \ge 0$. Mais, pour une hypersurface,
 ceci \'equivaut \`a demander que les courbures sectionnelles de $M$
 soient positives ou nulles \cite{Su}. 
Un th{\'e}or{\`e}me de Sacksteder \cite{S} assure que toute immersion
isom{\'e}trique de classe
$C^4$ d'une $3$-vari{\'e}t{\'e} {\`a} courbure sectionnelle 
non-n{\'e}gative est
forc{\'e}ment convexe et le r{\'e}sultat se ram{\`e}nerait au th{\'e}or{\`e}me
d'Archim{\`e}de. On a toutefois des r{\'e}sultats nouveaux en
affaiblissant la condition sur la courbure de Ricci. 

\medskip 

\noindent {\bf Proposition \ref{Ricci}.} {\it (n=3) 
Soit $\varphi~: M^3 \longrightarrow \RR^4$ une immersion isom{\'e}trique
de classe $C^2$ d'une vari{\'e}t{\'e} riemannienne lisse compacte sans bord.  
On suppose que $\varphi(M) \subset B^4(0,R)$. Si $M$ v{\'e}rifie $Ric
\ge - \alpha / R^2$ avec $ 0 < \alpha < 6$, alors 
$$\vol(M) <  \frac6{6-\alpha} R^3 T(M) $$
}

\demo  On utilise $\sum_i \langle \varphi, \, e_i \rangle ^2 
\le |\varphi|^2 \le R^2$ 
dans (\ref{n=3}) pour obtenir   
$$\text{vol}(M) \leq \int_M (p^3 - \frac23 pq^2) K dV + \frac\alpha{6}
\vol(M)$$
et de l{\`a} $\vol(M) \le 
\frac{6}{6-\alpha}R^3T(M)$ puisque $p^2+q^2=|\varphi|^2 \le
R^2$. L'{\'e}galit{\'e} n'est pas
atteinte puisque cela forcerait $q\equiv 0$ et $p\equiv R$. Dans ce
cas l'immersion devrait {\^e}tre le plongement standard de la sph{\`e}re
de rayon $R$, pour laquelle on a en fait
$\vol(\Ss^3(R))=R^3T(\Ss^3(R)) < \frac 6{6-\alpha}  R^3 T(\Ss^3(R))$.

\hfill{$\square$}

%%% Local Variables: 
%%% mode: latex
%%% TeX-master: "Gauss"
%%% TeX-master: "Gauss"
%%% End: 

%\setcounter{chapter}{3}
%\setcounter{section}{0}
%\setcounter{equation}{0}
\section{In{\'e}galit{\'e}s {\`a} caract{\`e}re local} \label{ineg locales}

Comme le
montre la section pr\'ec\'edente, l'in{\'e}galit{\'e}
(\ref{desir})  
est valable en petite dimension gr{\^a}ce {\`a} l'existence de formules
int{\'e}grales globales, qui reviennent en fait {\`a} la possibilit{\'e}
d'int{\'e}grer par parties~: une m\'ethode alternative 
pour obtenir les formules de Minkowski (\ref{minkowski}) est 
d'int{\'e}grer sur des ouverts constituant un d{\'e}coupage de
$M$ et de sommer. Les termes de bord s'annulent deux par deux \`a
cause du choix oppos\'e de l'orientation et, en dimension $2$, les
termes int\'erieurs sont controlables par la courbure totale. Dans
ce qui suit nous traitons les termes parasites qui apparaissent en
dimension $n\ge 3$ en simulant une int\'egration par parties o\`u la
contribution des termes de bord est compens\'ee par l'utilisation d'un 
\'epaisissement de l'ouvert relativement compact dont on estime le
volume. 

Tout ouvert sur lequel l'application de Gauss $G$ est
non-d\'eg\'en\'er\'ee et injective peut \^etre param\'etr\'e par
$G^{-1}$ qui, \`a son tour, s'exprime \`a l'aide de la
fonction support sur $\Ss^n$ et de sa
hessienne. La forme locale de l'in\'egalit\'e (\ref{desir}) s'exprimera
comme une in\'egalit\'e int\'egrale qu'on obtient par une technique
d'estimation emprunt\'ee \`a la th\'eorie des
op\'erateurs de type Monge-Amp\`ere \cite{Aubin, RT}. 
Dans la section \ref{a bas l'injectivite} nous renon{\c c}ons 
{\`a} l'hypoth{\`e}se d'injectivit{\'e} sur l'application de Gauss 
par un argument de recollement. La section \ref{local to global} 
utilise les formules int\'egrales d\'evelopp\'ees pr\'ec\'edememment 
pour obtenir (\ref{the best}) \`a partir de (\ref{the semi best}).
Nous obtenons aussi une
in{\'e}galit{\'e} {\`a} caract{\`e}re isop{\'e}rim{\'e}trique en dimension
deux. 

\subsection{Param{\'e}trisation d'Euler}\label{euler}
 
Consid{\'e}rons sur $M$ un domaine $U$ sur lequel l'application de Gauss  
$G:U\stackrel{\sim}{\longrightarrow} V\subset \Ss^n$ 
est injective. Son inverse  
$$\phi :V\stackrel{\sim}{\longrightarrow} U$$ 
est une param{\'e}trisation de $U$ et associe {\`a} un point $\theta$  
l'unique $x$ tel  
que le vecteur $\theta$ soit orthogonal {\`a} l'hyperplan $T_xM$ tangent  
{\`a} $M$ en $x$. Si $f(\theta)$ d{\'e}signe la distance orient{\'e}e entre 
l'origine $O$ et
$T_xM$, alors l'application $f$ est diff{\'e}rentiable sur $V$ et le  
gradient $\nabla f(\theta)\in T_\theta\Ss^n$ est bien d{\'e}fini. 
De plus, on aura   
$$\phi(\theta)=x=f(\theta)\cdot\theta+\nabla f(\theta)$$ 
apr{\`e}s des identifications naturelles avec des vecteurs de $\RR^{n+1}$. 
En g{\'e}om{\'e}trie convexe on appelle $f$ la  
``fonction support'' de $U$ et $\phi$ la  
``param{\'e}trisation d'Euler'' de l'hypersurface $U$. 
 
Soient $dV$ et $d\theta$ les {\'e}l{\'e}ments de volume sur $U\subset M$ et  
$V\subset\Ss^n$ respectivement. Alors 
$$dV=\pm\det(\phi_*(\theta))\cdot d\theta$$ 
On calcule maintenant  
$\phi_*:T_\theta\Ss^n\longrightarrow T_\theta\Ss^n\equiv
T_{\phi(\theta)}M^n\subset T_{\phi(\theta)}\RR^{n+1}$. On peut \'ecrire
$$\phi_*Y=Y(f)\theta+fY+D_Y(\nabla f)$$ 
Mais $T_{\phi(\theta)}M^n//T_\theta\Ss^n$ et en prenant ci-dessus la  
partie tangente {\`a} $\Ss^n$ on a 
$\phi_*Y=\phi_*^TY=fY+D_Y^T(\nabla f)=fY+D_Y^{\Ss^n}(\nabla f)$. 
Cela signifie pr{\'e}cis{\'e}ment que 
$$\phi_*=H_f+f\cdot Id$$ 
o{\`u} $H_f$ est la {\it hessienne} de $f$ sur la sph{\`e}re $\Ss^n$.  
C'est l'endomorphisme du fibr{\'e} tangent d{\'e}fini par 
$$H_f(\xi)=\nabla_\xi\nabla f$$ 
Ceci {\'e}tablit la formule (o{\`u} le signe d{\'e}pend des orientations) 
\begin{equation}\label{formule integrale}
\vol (U)=\pm\int_V\det(H_f+f\cdot Id)d\theta
\end{equation} 
Il est utile de voir l'application $\phi$ comme inverse de  
l'application de Gauss. L{\`a} o{\`u} cette derni{\`e}re est un  
diff{\'e}omorphisme, la d{\'e}riv{\'e}e $\phi_*$ sera l'inverse de
l'application   
de Weingarten et ses valeurs propres seront les inverses des courbures  
principales de $M$. En particulier la signature de $\phi_*$ sera la  
m{\^e}me que la signature de l'application de Weingarten. Sur les domaines  
de signature impaire, $\phi_*$ change l'orientation et la formule  
pr{\'e}c{\'e}dente s'{\'e}crit avec le signe ``$-$''. Sinon elle
s'\'ecrit avec le signe  
``$+$''.

\subsection{Coordonn{\'e}es st{\'e}r{\'e}ographiques} 
 
Dans cette section on exprime la formule int{\'e}grale pr{\'e}c{\'e}dente en  
coordonn{\'e}es st{\'e}r{\'e}ographiques. 
Soit $\Ss^n(\rho)=\{ x\in\RR^{n+1} :\vert x\vert =\rho\}$  
la sph{\`e}re de rayon $\rho$ dans $\RR^{n+1}$.  
Il est bien connu que la projection st{\'e}r{\'e}ographique du p{\^o}le  
nord induit sur $\RR^n$ une m{\'e}trique $g$ de composantes  
$$g_{ij}=\frac4{(1+\frac{\vert x\vert^2}{\rho^2})^2}\cdot\delta_{ij}$$ 
Les coefficients du tenseur de Christoffel sont 
$$\Gamma^i_{jk}=-\frac{2}{\rho^2+\vert x\vert^2} 
\{\delta_{ik}x^j+\delta_{ij}x^k-\delta_{jk}x^i\}$$ 
Pour faciliter l'{\'e}criture on note par $\partial_i$ les champs  
$\DP{}{x^i}$. Ce sont des champs orthogonaux dans la m{\'e}trique  
sph{\'e}rique.  
Le {\it gradient} sph{\'e}rique de $f$ satisfait  
$(\nabla f,\partial_i)=\DP{f}{x^i}$ et cela donne 
$$\nabla f_{\vert_x}=\frac14(1+\frac{\vert x\vert^2}{\rho^2})^2 
\cdot\DP{f}{x^i}\cdot\partial_i$$ 
La {\it hessienne} sph{\'e}rique de $f$ v{\'e}rifie 
$$\begin{array}{l} 
  H_f(\partial_i)\\ 
  =\nabla_{\partial_i}\nabla f=\nabla_{\partial_i}  
   \frac1{4\rho^4}(\rho^2+\vert x\vert^2)^2 
   \cdot\DP{f}{x^j}\cdot\partial_j\\ 
  =\frac1{4\rho^4}\cdot \{\partial_i(\rho^2+\vert x\vert^2)^2 
   \cdot\DP{f}{x^j}\cdot\partial_j  
   +(\rho^2+\vert x\vert^2)^2\cdot 
    \frac{\partial^2f}{\partial x^i\partial x^j}\cdot\partial_j  
   +(\rho^2+\vert x\vert^2)^2\cdot\partial_jf\cdot 
    \Gamma^l_{ij}\partial_l\}\\ 
  =\frac14(1+\frac{\vert x\vert^2}{\rho^2})^2\cdot\partial_{ij}f\cdot\partial_j 
   +\frac1{\rho^2}(1+\frac{\vert x\vert^2}{\rho^2})\cdot  
    x^i\partial_jf\cdot\partial_j 
    \\ \;\;\;\; 
   -\frac1{2\rho^2}(1+\frac{\vert x\vert^2}{\rho^2})\cdot\partial_jf\cdot 
    \{\delta_{li}x^j+\delta_{lj}x^i-\delta_{ij}x^l\}\partial_l\\ 
  =\frac14(1+\frac{\vert x\vert^2}{\rho^2})^2\cdot\partial_{ij}f\cdot\partial_j 
   +\frac1{2\rho^2}(1+\frac{\vert x\vert^2}{\rho^2})\cdot 
    [2x^i\partial_jf\cdot\partial_j-\partial_lf\{\delta_{ji}x^l+\delta_{lj}x^i- 
    \delta_{il}x^j\}\partial_j]\\ 
  =\frac14(1+\frac{\vert x\vert^2}{\rho^2})^2\cdot\partial_{ij}f\cdot\partial_j 
   +\frac1{2\rho^2}(1+\frac{\vert x\vert^2}{\rho^2})\cdot 
    [2x^i\partial_jf\cdot\partial_j- 
    \partial_lf\sum_{j\neq i} 
    \{\delta_{lj}x^i-\delta_{il}x^j\}\cdot\partial_j- 
    \partial_lf\cdot x^l\cdot\partial_i]\\ 
  =\frac14(1+\frac{\vert x\vert^2}{\rho^2})^2\cdot\partial_{ij}f\cdot\partial_j 
   +\frac1{2\rho^2}(1+\frac{\vert x\vert^2}{\rho^2})\cdot 
    [(x^i\partial_jf+x^j\partial_if)\cdot\partial_j- 
     (\nabla^{\RR^n} f_{\vert_x}\cdot x)\cdot\partial_i] 
\end{array} 
$$ 
Dans les formules qui vont suivre les gradients et hessiennes seront
implicitement pris par rapport {\`a} la  
m{\'e}trique euclidienne de $\RR^n$. Dans le cas contraire, on 
l'indiquera explicitement  
(par exemple $H_f^{\Ss^n}$ d{\'e}signera la hessienne de $f$ calcul{\'e}e
dans la m{\'e}trique sph{\'e}rique). 
On peut simplifier l'expression obtenue en posant 
\begin{equation}\label{g}
g=\frac{\rho^2+\vert x\vert^2}{2}\cdot f
\end{equation}
Alors  
$$\frac{\partial^2g}{\partial x^i\partial x^j}= 
 \frac{\rho^2+\vert x\vert^2}{2}\cdot \frac{\partial^2f}{\partial 
x^i\partial x^j}  
 +(x^i\partial_jf+x^j\partial_if)+ 
 f\cdot\delta_{ij}$$ 
et on obtient 
$$\begin{array}{l} 
 H_f^{\Ss^n}(x)(\partial_i)=
\frac1{2\rho ^2}(1+\frac{\vert x\vert ^2}{\rho ^2})\big[ \frac{\rho ^2
  + \vert x \vert ^2}{2} \partial_{ij} f \cdot \partial_j + (x^i
\partial_j f + x^j \partial_i f) \partial_j \big] - \frac1{2\rho ^2}
(1 + \frac{\vert x \vert ^2}{\rho ^2} ) (\nabla f\cdot x )\cdot
\partial_i \\
\quad \quad \quad \quad \quad
 =\frac1{2\rho^2}(1+\frac{\vert x\vert^2}{\rho^2})\cdot 
  \big[ H_g-(f+\nabla f\cdot x)\cdot Id \big] (\partial_i) 
 \end{array} 
$$ 
Soient maintenant $U$ et $V$ des domaines sur $M$ et respectivement $\Ss^n$ 
comme dans \S \ref{euler}. 
On choisit un point $N\in\Ss^n\setminus V$ et on note 
$\Omega$ l'image de $V$ par la projection st{\'e}r{\'e}ographique de $N$. 
Alors  
$$\int_V(H_f+f\cdot Id)d\theta= 
 \frac1{ 2^n \rho^{2n}} 
 \int_{\Omega}
   (1+\frac{\vert x\vert^2}{\rho^2})^n\cdot 
\det(H_g+(-f-\nabla f\cdot x+\frac{2\rho^2}{1+\frac{\vert 
x\vert^2}{\rho^2}}\cdot f))  
\cdot\frac{2^n}{(1+\frac{\vert x\vert^2}{\rho^2})^n}\;dx$$ 
donc 
\begin{center} 
\fbox{ $\displaystyle 
\vol (U)=\pm\frac1{\rho^{2n}}\int_{\Omega}\det(H_g+R_f\cdot Id)\;dx$}  
\end{center} 
o{\`u} on a not{\'e}  
$R_f=-f-\nabla f\cdot x+\frac{2\rho^4}{\rho^2+\vert x\vert^2}\cdot f$. 
Il est utile de remarquer que 
$$R_f=\frac2{\rho^2+\vert x\vert^2} \big[ (g-\nabla g\cdot x)+ 
  (\rho^4-\rho^2)\cdot f \big] $$ 
En particulier, pour $\rho=1$ on obtient 
\begin{equation}\label{Rg}
R_f=\frac2{1+\vert x\vert^2}\cdot(g-\nabla g\cdot x)
\stackrel{\text{not.}}{=} R_g
\end{equation}
et on pourra {\'e}crire 
 
\begin{center} 
\fbox{ $\displaystyle \vol (U)=\pm\int_{\Omega}\det(H_g+R_g\cdot Id)\;dx$} 
\end{center}

\subsection{Formes diff{\'e}rentielles} 
 
\noindent Consid{\'e}rons les 1-formes suivantes: 
$$\omega_i=d(\partial_ig)+R_g\cdot dx^i$$ 
On a $$\omega_1\wedge\ldots\wedge\omega_n= 
       \det(H_g+R_g\cdot Id)\;dx^1\wedge\ldots\wedge dx^n$$ 
donc 
$$\vol (U)=\pm\int_\Omega\omega_1\wedge\ldots\wedge\omega_n$$  
Rappelons que $R_g=\frac2{1+\vert x\vert^2}(g-\nabla g\cdot x)$  
et donc 
$$\begin{array}{l} 
  dR_g=-\frac4{(1+\vert x\vert^2)^2}(g-\nabla g\cdot x) 
   \cdot(x\cdot dx)-\frac2{1+\vert x\vert^2} 
   \cdot\sum_{i=1}^n x^i d(\partial_i g)\\ 
  \;\;\;\;\;\;\;  
  = - \frac2{1+\vert x\vert^2} \sum_{i=1}^n x^i\cdot\omega_i 
  \end{array} 
$$ 
Si on note $u:=\frac2{1+\vert x\vert^2}$ on a 
$$dR_g=-u\cdot \sum_{i=1}^n x^i\omega_i$$ 
et donc 
\begin{equation}\label{stability}
d\omega_i=dR_g\wedge dx^i=-u\sum_{j=1}^n x^j \omega_j \wedge dx^i
\end{equation}
C'est l'identit{\'e} (\ref{stability}) qui est la cl\'e de notre
r\'esultat. On pourra la penser comme
exprimant une sorte de stabilit{\'e} des $\omega_i$ par
diff{\'e}rentiation ext{\'e}rieure. 

\subsection{Une in{\'e}galit{\'e} locale} \label{Inegalites locales}
  
 Etablissons d'abord les notations. Pour des suites d'indices 
$I=(i_1,\ldots,i_k),\;J=(j_1,\ldots,j_{n-k})\subset\{1,\ldots ,n\}$ 
on note $\omega_I:=\omega_{i_1}\wedge\ldots\wedge\omega_{i_k}$ et 
$dx^J:=dx^{j_1}\wedge\ldots\wedge dx^{j_{n-k}}$. Tous les indices sont  
des entiers compris entre $1$ et $n$. 
 
\begin{prop}\label{locale}
  Soit $\Omega\subset\RR^n$ un ouvert connexe born{\'e} et $f:\Omega
  \longrightarrow \RR$ une fonction telle que la matrice  
$H_g+R_g\cdot Id$ soit non-d{\'e}g{\'e}n{\'e}r{\'e}e 
 sur $\Omega$, o{\`u} $g$ et $R_g$ sont d{\'e}finies par
  (\ref{g}) et (\ref{Rg}) avec $\rho$ {\'e}gal {\`a} $1$. Soient   
$I,J$ deux suites d'indices  
avec $\vert I\vert=k,\;\vert J\vert=n-k$. 
Pour tout compact $K\subset\Omega$ il existe une constante 
$C_k(K,\Omega)$ ind{\'e}pendante de $f$ telle que, pour toute fonction  
$\psi\in\mathcal{C}^\infty(\Omega)$ avec $supp(\psi)\subseteq K$, on ait: 
$$\vert\int_\Omega\psi\;\omega_I\wedge dx^J\vert\leq  
   C_k(K,\Omega)\cdot \vol(\Omega)\cdot\sup_K\vert\psi\vert\cdot 
      \sup_\Omega(\vert f\vert+\vert\nabla f\vert)^k$$ 
 
De plus, $C_k(K,\Omega)$ 
est proportionnelle par une constante ind{\'e}pendante de $K$ et
$\Omega$ {\`a} $(\frac{(1+r)^2}{\delta})^k$, o{\`u}  
$2 \delta=\text{\rm dist}(K,\partial\Omega)$ 
et $\Omega\subset\{\vert x\vert\leq r\}$. 
\end{prop} 
 
{\it D{\'e}monstration}:  Donnons d'abord deux majorations pour 
$\vert R_g\vert$ et $\vert\nabla g\vert$. 
 
Puisque $\rho=1$ on a  
$R_g=\frac{1-\vert x\vert^2}{1+\vert x\vert^2}\cdot f-\nabla f\cdot x$. 
Si $\Omega\subset B(0,r)=\{x\in\RR^n:\vert x\vert\leq r\}$ alors on trouve 
$$\vert R_g\vert\leq\vert f\vert + r\vert\nabla f\vert 
   \leq (1+r)\cdot(\vert f\vert+\vert\nabla f\vert)$$ 
 
De m{\^e}me, on a  
$\nabla g=f\cdot x+\frac{1+\vert x\vert^2}{2}\cdot\nabla f$ et si 
$\Omega\subset B(0,r)$ alors 
$$\vert\nabla g\vert\leq r\vert f\vert+\frac{1+r^2}{2}\vert\nabla f\vert 
   \leq \frac{(1+r)^2}{2}\cdot(\vert f\vert+\vert\nabla f\vert)$$ 
 
I) Traitons d'abord le cas o{\`u} $H_g+R_g\cdot Id$ est 
\underline{d{\'e}finie positive} sur $\Omega$. 
On proc{\`e}de par r{\'e}currence sur $k$. 
 Le premier pas de r{\'e}currence est $k=1$. On supposera  
$J=(2,\ldots,n)$ puisque les autres cas lui sont sym\'etriques. 
Deux situations se pr{\'e}sentent: $I\cap J=\emptyset$ et  
$I\cap J\neq\emptyset$. Les calculs suivent une id{\'e}e de  
Rauch et Taylor \cite{RT} reprise dans Aubin \cite{Aubin}. 
 
\begin{itemize} 
 \item{a)} $I\cap J=\emptyset$: on a successivement 
$$\begin{array}{l}\displaystyle 
\vert\int_\Omega \psi\;\omega_1\wedge dx^2\wedge\ldots\wedge dx^n\vert= 
\vert\int_\Omega \psi\;(\partial_{11}g+R_g)\; 
        dx^1\wedge\ldots\wedge dx^n\vert\\ 
\;\;\;\;\;\;\;\;\;\;\;\;\;\;\;\;\;\;\;\;\;\;\;\;\;\;\; 
  \leq \displaystyle 
       \sup_K\vert\psi\vert\int_K \vert(\partial_{11}g+R_g)\vert\; 
              dx^1\wedge\ldots\wedge dx^n\\ 
\;\;\;\;\;\;\;\;\;\;\;\;\;\;\;\;\;\;\;\;\;\;\;\;\;\;\;  
  =\displaystyle\sup_K\vert\psi\vert \int_K (\partial_{11}g+R_g)\; 
              dx^1\wedge\ldots\wedge dx^n  
  \end{array}$$ 
 
La derni{\`e}re {\'e}galit{\'e} fait usage de la positivit{\'e} de la matrice 
$H_g+R_g\cdot Id$. 
 
Soit maintenant $\gamma$ une fonction {\`a} support compact dans  
$\mathcal{C}^\infty(\Omega)$ telle que $\gamma_{\vert_K}\equiv 1$ et 
$0\leq\gamma\leq 1$. 
On peut la choisir telle que $\vert\nabla\gamma\vert\leq 1/\delta$. 
Alors 
$$\int_K(\partial_{11}g+R_g)\;dx\leq 
 \int_\Omega\gamma\;(\partial_{11}g+R_g)\;dx$$ 
 
$$\int_\Omega\gamma\partial_{11}g\;dx= 
-\int_\Omega\partial_1\gamma\cdot\partial_1g\;dx\leq  
\vol(\Omega)\cdot\sup_\Omega\vert\partial_1\gamma\vert\cdot 
  \sup_\Omega\vert\partial_1g\vert\leq 
\frac1\delta\cdot \vol(\Omega)\cdot \frac{(1+r)^2}{2} 
    \sup_\Omega(\vert f\vert+\vert\nabla f\vert)$$ 
 
$$\int_\Omega\gamma R_g\;dx\leq \vol(\Omega)\cdot(1+r)\cdot\sup_\Omega 
(\vert f\vert+\vert\nabla f\vert)$$ 
 
ce qui donne 
$$\vert\int_\Omega \psi\;\omega_1\wedge dx^2\wedge\ldots\wedge dx^n\vert 
\leq C(r,\frac1\delta)\cdot \vol(\Omega)\cdot 
  \sup_\Omega\vert\psi\vert\cdot  
  \sup_\Omega(\vert f\vert+\vert\nabla f\vert)$$ 
 
\item{b)} $I\cap J\neq\emptyset$: soit $i\in[2,n]$. 
$$\begin{array}{l}\displaystyle 
\vert\int_\Omega \psi\;\omega_i\wedge dx^2\wedge\ldots\wedge dx^n\vert= 
\vert\int_\Omega \psi\;\partial_{1i}g\; 
        dx^1\wedge\ldots\wedge dx^n\vert\\ 
\;\;\;\;\;\;\;\;\;\;\;\;\;\;\;\;\;\;\;\;\;\;\;\;\;\; 
  \leq \displaystyle\sup_K\vert\psi\vert\int_K \vert\partial_{1i}g\vert\; 
              dx^1\wedge\ldots\wedge dx^n\\ 
\;\;\;\;\;\;\;\;\;\;\;\;\;\;\;\;\;\;\;\;\;\;\;\;\;\;  
  \leq \displaystyle\frac12 \sup_K\vert\psi\vert   
       \int_K ((\partial_{11}g+R_g)+(\partial_{ii}g+R_g))\; 
              dx^1\wedge\ldots\wedge dx^n\\ 
\;\;\;\;\;\;\;\;\;\;\;\;\;\;\;\;\;\;\;\;\;\;\;\;\;\;  
  \leq \displaystyle C(r,\frac1\delta)\cdot \vol(\Omega)\cdot 
  \sup_\Omega\vert\psi\vert\cdot  
  \sup_\Omega(\vert f\vert+\vert\nabla f\vert) 
  \end{array}$$ 
\end{itemize}  
La derni{\`e}re in{\'e}galit{\'e} est une cons{\'e}quence
imm{\'e}diate du point a),  
tandis que le fait crucial 
$$\vert\partial_{1i}g\vert\leq 
    \frac12((\partial_{11}g+R_g)+(\partial_{ii}g+R_g))$$ 
d{\'e}coule de la positivit{\'e} de la matrice $H_g+R_g\cdot Id$. 
Cette derni{\`e}re in{\'e}galit{\'e} reste valable pour des mineurs  
d'ordre quelconque  
d'une matrice d{\'e}finie positive. Consid{\'e}rons deux suites d'indices 
$I=(i_1,\ldots,i_k),\;J=(j_1,\ldots,j_{n-k})\subset\{1,\ldots ,n\}$.  
On note $\stackrel{\sim}{J}$ la suite compl{\'e}mentaire de $J$.  
Pour une matrice $M=(m_{ij})\in\mathcal{M}_n(\RR)$ on note  
$M[I,\stackrel{\sim}{J}]$  
la matrice carr{\'e}e  
construite sur les lignes $I$ et les colonnes $\stackrel{\sim}{J}$ de 
$M$ et $M_{I,\stackrel{\sim}{J}}$ son d{\'e}terminant (par exemple, pour $k=1$ on a  
$M_{ij}=m_{ij}$). 
 
Si $M$ est d{\'e}finie positive alors  
 
$$\forall \; i,j \in\overline{1,n}\;:\;\vert M_{ij}\vert\leq 
\frac12(M_{ii}+M_{jj})$$ 
puisque $M_{\{i,j\},\{i,j\}}\geq 0$. De plus, une fois l'orientation  
fix{\'e}e sur $\RR^n$ (repr{\'e}sent{\'e}e - par exemple - par la base  
canonique 
$(e_1,\ldots,e_n)$), $M$ d\'etermine naturellement un endomorphisme  
{\it d{\'e}fini positif} sur $\Lambda^k(\RR^n)$ muni de la m{\'e}trique  
induite par $\RR^n$. Sa matrice dans la base  
$e_{i_1}\wedge\ldots\wedge e_{i_k}$ a comme entr{\'e}es justement les  
$M_{I,\stackrel{\sim}{J}}$. On trouve donc l'in{\'e}galit{\'e} 
$$\forall\;I,J,\vert I\vert=k,\vert J\vert=n-k\;:\; 
\vert
M_{I,\stackrel{\sim}{J}}\vert\leq\frac12(M_{I,I}+
M_{\stackrel{\sim}{J},\stackrel{\sim}{J}})$$ 
 
On suppose maintenant l'in{\'e}galit{\'e} vraie pour $k-1$ et on la prouve  
pour $k$. Les hypoth{\`e}ses sont $supp\;\psi\subseteq K$, $\vert I\vert=k$, 
$\vert J\vert=n-k$. On note $M:=H_g+R_g\cdot Id$. Alors  
 
$$\begin{array}{l}\displaystyle 
\vert\int_\Omega \psi\;\omega_I\wedge dx^J\vert= 
\vert\int_\Omega \psi\;M_{I,\stackrel{\sim}{J}}\; 
        dx^1\wedge\ldots\wedge dx^n\vert\\ 
\;\;\;\;\;\;\;\;\;\;\;\;\;\;\;\;\;\;\;\;\;\;\;\;\;\; 
  \leq \displaystyle\sup_K\vert\psi\vert\int_K \vert M_{I,\stackrel{\sim}{J}}\vert\; 
              dx^1\wedge\ldots\wedge dx^n\\ 
\;\;\;\;\;\;\;\;\;\;\;\;\;\;\;\;\;\;\;\;\;\;\;\;\;\;  
  \leq\displaystyle\frac12 \sup_K\vert\psi\vert   
       \int_K (M_{I,I}+M_{\stackrel{\sim}{J},\stackrel{\sim}{J}})\; 
              dx^1\wedge\ldots\wedge dx^n\\ 
\;\;\;\;\;\;\;\;\;\;\;\;\;\;\;\;\;\;\;\;\;\;\;\;\;\;  
  =\displaystyle\frac12\sup_K\vert\psi\vert 
  (\int_K\omega_I\wedge dx^{\stackrel{\sim}{I}}+\int_K\omega_{\stackrel{\sim}{J}}\wedge dx^J) 
  \end{array}$$ 
 
Il suffit donc de trouver les estimations pour $J=\stackrel{\sim}{I}$. Comme avant,  
soit $\gamma$ une fonction {\`a} support compact dans  
$\mathcal{C}^\infty(\Omega)$ telle que $\gamma_{\vert_K}\equiv 1$, 
$0\leq\gamma\leq 1$ et $dist(supp\;\gamma,\partial\Omega)\geq\delta/2$. 
On peut la choisir telle que $\vert\nabla\gamma\vert\leq 1/\delta$. Alors 
 
$$\begin{array}{l}\displaystyle  
   \int_K\omega_I\wedge dx^{\stackrel{\sim}{I}}\leq 
    \int_\Omega\gamma\;\omega_I\wedge dx^{\stackrel{\sim}{I}}\\ 
\;\;\;\;\;\;\;\;\;\;\;\;\;\;\;\;\;\;\; 
   =\displaystyle 
     \int_\Omega\gamma\;d(\partial_{i_1}g)\wedge\omega_{I-\{i_1\}}\wedge 
     dx^{\stackrel{\sim}{I}}+ 
    \int_\Omega\gamma\;R_g\;dx^{i_1}\wedge\omega_{I-\{i_1\}}\wedge 
     dx^{\stackrel{\sim}{I}} 
\end{array}$$ 
 
Par l'hypoth{\`e}se de r{\'e}currence on a  
$$\begin{array}{l}\displaystyle 
  \int_\Omega\gamma\;R_g\;dx^{i_1}\wedge\omega_{I-\{i_1\}}\wedge 
     dx^{\stackrel{\sim}{I}}\leq  
    C_{k-1}(r,\frac4\delta)\cdot\sup_\Omega\vert\gamma\cdot R_g\vert\cdot 
    \vol(\Omega)\cdot 
    \sup_\Omega(\vert f\vert+\vert\nabla f\vert)^{k-1}\\ 
\;\;\;\;\;\;\;\;\;\;\;\;\;\;\;\;\;\;\;\;\;\;\;\;\;\;\;\;\;\;\;\; 
\;\;\;\;\;\;\;\;\;\;\; 
  \leq \displaystyle(1+r)\cdot 4^{k-1}\cdot C_{k-1}(r,\frac1\delta)\cdot 
         \vol(\Omega)\cdot 
         \sup_\Omega(\vert f\vert+\vert\nabla f\vert)^k 
  \end{array} 
$$ 
 
Apr{\`e}s une int{\'e}gration par parties on obtient aussi 
$$\begin{array}{l}\displaystyle 
  \int_\Omega\gamma\;d(\partial_{i_1}g)\wedge\omega_{I-\{i_1\}}\wedge 
     dx^{\stackrel{\sim}{I}}\\ 
   \displaystyle 
   =-\int_\Omega\partial_{i_1}g\;d\gamma\wedge\omega_{I-\{i_1\}}\wedge 
     dx^{\stackrel{\sim}{I}} 
     -\sum_{l=2}^k\int_\Omega\gamma\;\partial_{i_1}g\; 
      \omega_{i_2}\wedge\ldots\wedge dR_g\wedge dx^{i_l}\wedge\ldots 
      \wedge\omega_{i_k}\wedge dx^{\stackrel{\sim}{I}}\\ 
   \displaystyle 
   =-\int_\Omega\partial_{i_1}g\;(\partial_s\gamma\;dx^s) 
    \wedge\omega_{I-\{i_1\}}\wedge dx^{\stackrel{\sim}{I}} 
    +\sum_{l=2}^k\int_\Omega\gamma\;\partial_{i_1}g\; 
                                 \frac2{1+\vert x\vert^2}\; 
      \omega_{i_2}\ldots(x^s\;\omega_s)\wedge dx^{i_l} 
      \ldots\omega_{i_k}\wedge dx^{\stackrel{\sim}{I}}\\ 
   \displaystyle 
   \leq nC_{k-1}(r,\frac4\delta)\;\frac{(1+r)^2}2\; 
     \sup_\Omega(\vert f\vert+\vert\nabla f\vert)\;\frac1\delta\;  
    \vol(\Omega)\sup_\Omega(\vert f\vert+\vert\nabla f\vert)^{k-1}\\ 
    \displaystyle 
\;\;\;\;\;\;\;\;\;\;\;\;\;\;\;\; 
    +(k-1)C_{k-1}(r,\frac4\delta)\;\frac{(1+r)^2}2\; 
     \sup_\Omega(\vert f\vert+\vert\nabla f\vert)\; 
   \vol(\Omega)\sup_\Omega(\vert f\vert+\vert\nabla f\vert)^{k-1}\\ 
    \displaystyle 
\leq (n+k-1)\;\frac{(1+r)^2}2\;\frac1\delta\;4^{k-1} 
         \;C_{k-1}(r,\frac1\delta)\; 
      vol_{\RR^n}(\Omega)\sup_\Omega(\vert f\vert+\vert\nabla f\vert)^k 
  \end{array} 
$$ 
 
Cela donne l'estimation finale 
$$\int_K\omega_I\wedge dx^{\stackrel{\sim}{I}}\leq 
    ((n+k-1)\;\frac{(1+r)^2}2+1+r)\;4^{k-1}\;\frac1\delta 
         \;C_{k-1}(r,\frac1\delta)\; 
    \vol(\Omega)\sup_\Omega(\vert f\vert+\vert\nabla f\vert)^k$$ 
 
ou bien 
$$\vert\int_K\psi\;\omega_I\wedge dx^{\stackrel{\sim}{I}}\vert\leq 
   C_k(r,\frac1\delta)\;\sup_\Omega\vert\psi\vert\; 
   \vol(\Omega)\sup_\Omega(\vert f\vert+\vert\nabla f\vert)^k$$ 
 
avec la constante {\'e}quivalente {\`a} $(\frac1\delta)^k$ et $r^{2k}$ 
$$C_k(r,\frac1\delta)=((n+k-1)\;\frac{(1+r)^2}2+1+r)\;4^{k-1}\;\frac1\delta 
         \;C_{k-1}(r,\frac1\delta)$$ 
 
Cela ach{\`e}ve la d{\'e}monstration du cas d{\'e}fini positif. 
 
II) On traite maintenant le cas d'une signature arbitraire. 
 
{\bf Fait} simple et fondamental: Soit $\Omega$ un ouvert  
connexe et $M(x)$ une famille  
de matrices sym{\'e}triques d{\'e}pendant contin{\^u}ment du  
param{\`e}tre $x\in\Omega$. On a l'{\'e}quivalence: 
 
i) $M(x)$ est d{\'e}finie positive pour tout $x\in \Omega$ 
 
ii) Pour tout $x\in\Omega$ on a $\det M(x)\neq 0$ et il existe  
un $x_0\in\Omega$ tel que $M(x_0)$ est d{\'e}finie positive. 
 
Appliquons cela dans notre cas: choisissons un point $x_0\in\Omega$.  
Il existe une matrice diagonale 
$A=\text{\rm diag}(\epsilon_1,\ldots,\epsilon_n)$  
telle que $\epsilon_i=\pm1$ et 
$A\cdot(H_g+R_g\cdot Id)_{\vert_{x_0}}>0$. Vu que $H_g+R_g\cdot Id$ 
est non-d{\'e}g{\'e}n{\'e}r{\'e}e sur $\Omega$ on d{\'e}duit par le fait  
pr{\'e}c{\'e}dent que $A\cdot (H_g+R_g\cdot Id)>0$ sur $\Omega$. 
 
Multiplier par $A$ signifie multiplier certaines lignes de  
$H_g+R_g\cdot Id$ par $\pm1$. Vu que tous les d{\'e}terminants ont {\'e}t{\'e} 
exprim{\'e}s par les formes $\omega_i$ qui correspondent aux lignes de  
$H_g+R_g\cdot Id$, cela revient {\`a} multiplier certaines des formes  
$\omega_i$ par $\pm1$. Le cas d{\'e}fini positif nous donne l'in{\'e}galit{\'e} 
$$\vert\int_\Omega\psi\;\stackrel{\sim}\omega_I\wedge dx^J\vert\leq  
   C_k(K,\Omega)\cdot \vol(\Omega)\cdot\sup_K\vert\psi\vert\cdot 
      \sup_\Omega(\vert f\vert+\vert\nabla f\vert)^k$$ 
o{\`u}  
$\stackrel{\sim}\omega_I=\stackrel{\sim}\omega_{i_1}\wedge\ldots\wedge\stackrel{\sim}\omega_{i_k}$ 
et $\stackrel{\sim}\omega_i=\epsilon_i\cdot\omega_i=\pm\omega_i$. Mais on a de  
fa{\c c}on {\'e}vidente 
$$\int_\Omega\psi\;\omega_I\wedge dx^J= 
    \pm\int_\Omega\psi\;\stackrel{\sim}\omega_I\wedge dx^J$$ 
ce qui donne l'in{\'e}galit{\'e} d{\'e}sir{\'e}e pour une signature quelconque  
et ach{\`e}ve la d{\'e}monstration. 
 
\hfill{$\square$} 

\noindent La proposition pr\'ec\'edente a \'et\'e d\'emontr\'ee en vue du 
 
\begin{cor} 
  Avec les notations de la proposition \ref{locale}, si 
$H_g+R_g\cdot Id$ est non-d{\'e}g{\'e}n{\'e}r{\'e}e 
sur l'ouvert born{\'e} $\Omega\subset\{\vert x\vert\leq r\}$ alors 
pour tout compact $K\subseteq\Omega$  
il existe une constante $C(K,\Omega)$ telle que  
$$\int_K\omega_1\wedge\ldots\wedge\omega_n\leq 
    C(K,\Omega)\cdot \vol(\Omega)\cdot 
    \sup_\Omega(\vert f\vert+\vert\nabla f\vert)^n$$ 
Cette constante est proportionnelle {\`a}  
$(\frac{(1+r)^2}{\delta})^n$, o{\`u} 
$\delta=\text{\rm dist}(K,\partial\Omega)$  
\end{cor} 
 
Dans l'in{\'e}galit{\'e} ci-dessus, le terme de gauche repr{\'e}sente le
volume d'un compact de $M$ sur lequel l'application de Gauss est
non-d{\'e}g{\'e}n{\'e}r{\'e}e et injective, tandis que le membre de droite
repr{\'e}sente le volume d'un {\'e}paississement par $\delta$ de son image
dans $\RR^n$ par la compos{\'e}e de la projection
st{\'e}r{\'e}ographique avec l'application de Gauss. 
Si $r\le 1$, la distance euclidienne est comparable {\`a} la distance
mesur{\'e}e sur $\Ss^n$ et on en d{\'e}duit le

\begin{cor}\label{coro}
  Soit $\varphi:M \longrightarrow E^{n+1}$ une immersion isom\'etrique d'hypersurface dont l'image est contenue dans la boule $B^{n+1}(0,R)$. 
Soit $U\subset M$ un ouvert sur lequel l'application de Gauss est
  non-d{\'e}g{\'e}n{\'e}r{\'e}e, injective et d'image contenue dans une
  demi-sph{\`e}re. Il existe une constante $C_n$ ind{\'e}pendante de $U$
  et de $M$ telle que, si $K \subset U$ est un compact, on a
 \begin{equation}\label{locale intrinseque}
  \vol \: (K) \le C_n R^n\frac1{\big(\text{d}_{\: \Ss^n}(K,\, \partial
    U)\big)^n}  \, T(U) 
 \end{equation}
 Ici $\text{d}_{\Ss^n}(K,\, \partial U)$ repr{\'e}sente la distance
 entre $K$ et $\partial U$ mesur{\'e}e sur $\Ss^n$ via l'application de
 Gauss, c'est-{\`a}-dire
 \begin{equation}\label{distance spherique}
 d_{\Ss^n}(K,\partial U) \ \stackrel{\text{\rm d{\'e}f.}}{=}  \ 
   \text{\rm dist}_{\Ss^n}\big(Gauss(K),\ Gauss(\partial U)\big) \ = \
\min_{x\in K} \ 
 \text{\rm dist}_{\Ss^n}\big(Gauss(x),\ Gauss(\partial U)\big)
 \end{equation}
\end{cor}

\subsection{Globalisation}\label{a bas l'injectivite}

Dans ce qui suit nous expliquons comment on peut {\'e}liminer les deux 
derni{\`e}res hypoth{\`e}ses sur l'application de Gauss dans le corollaire
\ref{coro}. L'application de Gauss sera not\'ee $G$
ou $Gauss$ et la courbure totale $T(U)$ d'un domaine $U$ sera parfois 
not\'ee $\vol(Gauss_U)$. Donnons d'abord la

\begin{defi} \label{la distance locale} 
Soit $U\subset M$ un ouvert sur lequel l'application de
  Gauss est non-d{\'e}g{\'e}n{\'e}r{\'e}e et $K\subset U$ un compact. La
  distance {\rm sph{\'e}rique locale} entre $K$ et $\partial U$ est
  d{\'e}finie comme 
 \begin{eqnarray}\label{distance locale}
 \bar{d}_{\: \Ss^n}(K,\, \partial U) & \stackrel{\text{\rm
 d{\'e}f.}}{=} &
 \min_{x\in K} \ 
 \sup_{
        \begin{array}{c} 
          U \supseteq \mathcal{U}\ni x \\ 
          Gauss(\mathcal{U}) \subset \text{\rm \small demi-sph{\`e}re} \\
          Gauss_{\vert _ \mc{U}} \ \text{\rm \small injective}
        \end{array}
       } 
           \text{d}_{\: \Ss^n}\big(x,\, \partial \mathcal{U}\big) \\
 & = & 
\min_{x\in K} \ 
 \sup_{
        \begin{array}{c} 
          U \supseteq \mathcal{U}\ni x \\ 
          Gauss(\mathcal{U}) \subset \text{\rm \small demi-sph{\`e}re} 
        \end{array}
       } 
           \text{d}_{\: \Ss^n}\big(x,\, \partial \mathcal{U}\big)
\nonumber 
 \end{eqnarray}

\end{defi}

\noindent Les deux quantit\'es d\'efinies ci-dessus sont \'egales en
vues de la non-d\'eg\'en\'erescence de $Gauss$ sur $U$. Pour montrer
l'\'egalit\'e il suffit de voir que, pour $x\in K$ fix\'e et $\mc{U}$
satisfaisant 
$x\in \mc{U} \subset U$ et $Gauss(\mc{U}) \subset \text{\rm
demi-sph\`ere}$, il existe un $x\in \mc{U}'\subset \mc{U}$ sur lequel
$Gauss$ est injective et tel que $d_{\Ss^n}(x, \
\partial\mc{U})=d_{\Ss^n}(x, \ \partial\mc{U}')$. Or, par d\'efinition 
de $d_{\Ss^n}(x, \ \partial \mc{U})$, l'application
$G_{\vert_{\mc{U}}}$ est propre au dessus de la boule ouverte
$B_{\Ss^n}\big(G(x), \ d_{\Ss^n}(x, \ \partial \mc{U}) \big)$. Comme
c'est aussi un diff\'eomorphisme local on d\'eduit que c'est un
rev\^etement. L'ouvert $\mc{U}'$ cherch\'e sera la feuille qui
contient $x$. 

\begin{thm}\label{theoreme local principal}
  Soit $\varphi:M \longrightarrow E^{n+1}$ une immersion isom\'etrique
d'hypersurface dont l'image est contenue dans la boule $B^{n+1}(0,R)$.
Soit $U\subset M$ un ouvert sur lequel l'application de Gauss est
  non-d{\'e}g{\'e}n{\'e}r{\'e}e. Il existe une constante $\bar{C}_n$ 
  ind{\'e}pendante de $U$
  et de $M$ telle que, si $K \subset U$ est un compact, on a
  \begin{equation}\label{gener}
 \vol \: (K) \le \bar{C}_n R^n \frac1{\big(\bar{d}_{\: \Ss^n}(K,\,
 \partial U)\big)^n} 
 \, \vol \: (Gauss_U) 
 \end{equation}
 Ici $\vol \: (Gauss_U)$ repr{\`e}sente le volume recouvert sur
 $\Ss^n$ par
 l'application de Gauss restreinte {\`a} $U$ 
 {\rm compt{\'e} avec multiplicit{\'e}s}.
\end{thm}

Avant de prouver le th{\'e}or{\`e}me nous donnons quelques
propri{\'e}t{\'e}s de la distance $\bar{d}_{\: \Ss^n}(K,\, \partial U) $.

\renewcommand{\theenumi}{\roman{enumi}}

\begin{enumerate}
\item si l'application de Gauss est non seulement
  non-d{\'e}g{\'e}n{\'e}r{\'e}e sur $U$, mais aussi injective et d'image
  contenue dans une demi-sph{\`e}re (i.e. sous les hypoth{\`e}ses du
  corollaire \ref{coro}), on a $\bar{d}_{\: \Ss^n}(K,\, \partial U)=
  d_{\: \Ss^n}(K,\, \partial U) $ ;

\item on a toujours $\bar{d}_{\: \Ss^n}(K,\, \partial U)  > 0$ et
  $$d_{\: \Ss^n}(K,\, \partial U) \le
  \bar{d}_{\: \Ss^n}(K,\, \partial U) \le \frac\pi2 $$
  Il se peut toutefois que $d_{\: \Ss^n}(K,\,
  \partial U)$ soit nulle si on n'impose pas l'injectivit\'e de $G$
sur $U$, \`a savoir lorsque l'image de $\partial U$ recoupe l'image de 
$K$. 

\item 
$$\bar{d}_{\: \Ss^n}\big(K,\, \ \partial U \big) \ = \ \bar{d}_{\:
    \Ss^n}\big(\partial K,\, \ \partial U \big)$$

\item d{\'e}finissons une semidistance sur $M$ par 
    $$ d_{\Ss^n}(x,y) \ \stackrel{\text{\rm d{\'e}f.}}{=} \ \text{\rm
      dist}_{\Ss^n} \big( Gauss(x), \ Gauss(y) \big) $$
  Cette semi-distance est non-d{\'e}g{\'e}n{\'e}r{\'e}e sur un ouvert $U$ si
  et seulement si l'application de Gauss est injective sur $U$. 
  La notation $d_{\Ss^n}$ est justifi\'ee par l'identit\'e 
    $$d_{\Ss^n} \big( K, \ \partial U \big) \ = \ \min_{x\in K,\ y\in
      \partial U} d_{\Ss^n}(x,y)$$

  \noindent 
  Pour une application de Gauss non-injective les boules 
  $$B_{d_{\Ss^n}}(x,r) \ = \ \big\{ y\in M \ : \ d_{\Ss^n}(x,y)<r
  \big\}, \qquad r>0$$
  \noindent sont en g{\'e}n{\'e}ral disconnexes. D{\'e}finissons la {\it boule
    connexe} de rayon $r$ centr{\'e}e en $x$ comme 
  $$\til{B}_{d_{\Ss^n}}(x,r) \ \stackrel{\text{\rm d{\'e}f.}}{=} \
  \text{\rm la composante connexe de $x$ dans $B_{d_{\Ss^n}}(x,r)$},
  \qquad r>0$$
  avec la convention   $\til{B}_{d_{\Ss^n}}(x,0)=\{x\}$.
  D{\'e}finissons aussi le {\it rayon sph\'erique 
  d'injectivit{\'e} de l'application
  de Gauss} en $x$ comme 
  $$R_x \ \stackrel{\text{\rm d{\'e}f.}}{=} \
   \sup \big\{ r \ : \ d_{\Ss^n} \text{ \rm non-d{\'e}g{\'e}n{\'e}r{\'e}e sur
   $\til{B}_{d_{\Ss^n}}(x,r) $} \big\}$$
 Si l'application de Gauss est non-d{\'e}g{\'e}n{\'e}r{\'e}e en $x$ alors
 $R_x>0$ par injectivit\'e locale. Le lemme qui suit sera fondamental
dans la preuve du th{\'e}or{\`e}me \ref{theoreme local principal}.
\end{enumerate}

 \begin{lem}\label{lemme}
  Soit $U$ un ouvert sur lequel l'application de Gauss est
  non-d{\'e}g{\'e}n{\'e}r{\'e}e et $K \subset U$ un compact. Pour tout
  $x\in K$ on a $R_x \ge \bar{d}_{\Ss^n} \big( K, \ \partial U \big)$.
 \end{lem}
  
 {\small \it D{\'e}monstration du lemme:} L'{\'e}nonc{\'e} du lemme
 est en fait une reformulation g{\'e}om{\'e}trique des d{\'e}finitions. 
 Soit $x \in K$. Il suffit de
 montrer  
  $$R_x \ge \sup _{
        \begin{array}{c} 
          U \supseteq \mathcal{U}\ni x \\ 
          Gauss_{\vert _ \mc{U}} \ \text{\rm \small injective}
        \end{array}
       } 
           d_{\: \Ss^n} \big(x,\, \partial \mathcal{U}\big)$$ 
  Soit $\mc{U}$ tel que $Gauss_{\vert_{\mc{U}}}$ injective. Alors, par
           d{\'e}finition,  
  $\til{B}_{d_{\Ss^n}}\big(x, \, d_{\Ss^n}(x,\ \partial \mc{U}) \big)
           \subseteq \mc{U}$ et donc $R_x \ge d_{\Ss^n}(x,\ \partial
           \mc{U})$, ce qui ach{\`e}ve la preuve. 

\hfill{$\square$}
  
De ce fait on peut regarder la quantit{\'e} $\bar{d}_{\Ss^n}
           \big( K, \ \partial U \big)$ comme fournissant une borne
           inf{\'e}rieure uniforme sur $K$ pour 
           les rayons d'injectivit{\'e} de
           l'application de Gauss. 

\medskip 

\noindent
{\small \it D{\'e}monstration du th{\'e}or{\`e}me \ref{theoreme local
    principal} :} 
Il suffit de faire la d\'emonstration pour $R=1$. 
Nous allons fournir un argument de recollement. Il est possible
de sous-diviser $K$ en compacts $K'$ et de choisir autour de chaque
$K'$ un ouvert $U'\subset U$ satisfaisant les propri{\'e}t{\'e}s suivantes:

- la restriction de 
  l'application de Gauss {\`a} chaque $U'$ est injective et d'image
  contenue dans une demi-sph{\`e}re;

- $  d_{\: \Ss^n}(K',\,
  \partial U') \ge \frac1C \bar{d}_{\: \Ss^n}(K,\, \partial U) $, avec
  $C$ une constante universelle ind{\'e}pendante de la dimension;

- une intersection non-vide d'ouverts $U'$ compte au plus $D_n$
  {\'e}l{\'e}ments, avec $D_n$ une constante d{\'e}pendant
  uniquement de la dimension.

Une fois une telle sous-division de $K$ construite, on applique 
le corrolaire \ref{coro} pour chaque couple $(K',\, U')$ et on trouve:
\begin{eqnarray*}
 \vol(K)=\sum \vol(K') & \le & C_n \sum \frac1{\big( d_{\: \Ss^n}(K',\,
  \partial U') \big)^n} \vol(Gauss_{U'}) \\
 & \le & C_n  \frac{C^n}{\big(  \bar{d}_{\: \Ss^n}(K,\, \partial U)
  \big)^n } \sum \vol(Gauss_{U'}) \\
 & \le & D_n C_n \frac{C^n}{\big(  \bar{d}_{\: \Ss^n}(K,\, \partial U)
  \big)^n } \vol(Gauss_U) \\
 & = & \bar{C}_n \frac{1}{\big(  \bar{d}_{\: \Ss^n}(K,\, \partial U)
  \big)^n } \vol(Gauss_U)
\end{eqnarray*}
avec $\bar{C}_n=C^n D_n C_n$. 

Remarquons d'abord qu'il suffit de faire la d{\'e}monstration dans le
cas o{\`u} $Gauss$ est non-d{\'e}g{\'e}n{\'e}r{\'e}e sur $\overline{U}$: on
pourra ensuite consid{\'e}rer une exhaustion de $U$ par des ouverts
relativements compacts contenant $K$. Par continuit{\'e} de la distance
$\bar{d}_{\: \Ss^n}$ l'in{\'e}galit{\'e} pour $U$ lui-m{\^e}me s'ensuivra. 
On travaillera par la suite sous cette hypoth{\`e}se. 

La construction de la sous-division
$\{(K', \, U') \}$ se fait en deux {\'e}tapes:

{\underline{Etape 1:}} r{\'e}duction au cas o{\`u} $Gauss(U)$ est contenu
dans une demi-sph{\`e}re.  Fixons une triangulation $\{ T_i \}$ de la
sph{\`e}re $\Ss^n$ telle que chaque simplexe soit contenu dans une boule
de rayon $\frac{3\pi}{8}$. Pour chaque $T_i$ consid{\'e}rons le
voisinage $V_i = \{ p\in\Ss^n \, : \, \text{\rm dist}_{\Ss^n}(p,\, 
T_i)<\frac\pi8 \}$. Tout $V_i$ sera contenu dans une demi-sph{\`e}re et
chaque $V_i$ 
intersectera  au plus $c_n$ autres $V_j$, avec $c_n$ une constante qui
d{\'e}pend de la triangulation. La constante $c_n$
sera d{\'e}sormais fix{\'e}e, et il est facile d'obtenir des estimations
sur $c_n$ en fonction de $n$ en construisant une triangulation
explicite. 

On construit maintenant une division de $K$ et $U$ en prenant des
pr{\'e}images par $Gauss$. 
Soient $U_i^l$ les composantes connexes de 
$G^{-1}(V_i) \cap U$. Puisque $G$ est un diff{\'e}omorphisme local au
voisinage de  
$\overline{U}$, celles-ci co{\"\i}ncident avec les int{\'e}rieurs 
des composantes connexes de $G^{-1}(\overline{V_i}) \cap
\overline{U}$, et ces derni{\`e}res sont en nombre fini par compacit{\'e}
de $\overline{U}$. Donc les $U_i^l$ sont en nombre fini.  

Soient $K_i^l=G^{-1}(T_i) \cap K \cap U_i^l$. Les $K_i^l$ ne sont pas
n{\'e}cessairement connexes, mais ils sont d'int\'erieurs disjoints et on a 
$$\bigcup_l K_i^l = G^{-1}(T_i) \cap K$$

Ainsi $\{ (K_i^l,U_i^l) \}_{i,l}$ forme une sous-division de $(K,U)$
et $G(U_i^l)$ est inclus dans une demi-sph{\`e}re par
construction. De plus 
$$\bar{d}_{\Ss^n}(K_i^l, \, \partial U_i^l)= \min \{
\bar{d}_{\Ss^n}(K_i^l, \, \partial U), \, \frac\pi8 \} \ge 
\frac14 \bar{d}_{\Ss^n}(K, \, \partial U)$$
et le nombre de $U_i^l$ qui s'intersectent {\`a} la fois est born{\'e} par
$c_n$. 

Il suffira ainsi de trouver pour chaque couple $(K_i^l, \, U_i^l)$ une
sous-division avec les propri{\'e}t{\'e}s d{\'e}sir{\'e}es. Si les constantes
qui appara{\^\i}tront sont not{\'e}es $C'$ et $D'_n$, alors la
construction aura \'et\'e faite dans le cas g{\'e}n{\'e}ral avec $C=4C'$,
$D_n=c_n\cdot D'_n$. 

{\it Remarque:} Pour minimiser la constante finale il faudra trouver
un bon rapport entre le $c_n$ et la largeur de l'{\'e}paisissement $V_i$
de chaque $T_i$. La constante finale sera meilleure avec un $c_n$
petit et un {\'e}paisissement large. Or ces deux demandes sont
antagoniques et le produit des deux quantit{\'e}s  
est minimis{\'e} par un certain choix de
triangulation. Dans notre cas, on ne s'int{\'e}resse pas {\`a} la
meilleure constante et on cherche juste une information qualitative.

{\underline{Etape 2:}} on suppose d{\'e}sormais que $Gauss(U)$ est
inclus dans une demi-sph{\`e}re. Soit 
$$\delta =\frac1{15} \ \bar{d}_{\Ss^n}(K, \, \partial U)$$
Le lemme \ref{lemme} assure que $Gauss$ est injective sur 
$\til{B}(x, \, 15 \delta) \stackrel{\text{\rm not.}}{=} \
\til{B}_{d_{\Ss^n}}(x, \, 15 \delta)$ pour tout $x\in K$. 
Soit $\big( \til{B}(x_i, \, 3 \delta) \big)_{i=1, N}$ un
recouvrement fini {\it minimal} de $K$ par des boules de rayon $3
\delta$ centr{\'e}es en des points de $K$. On pose 
$$B_i=\til{B}(x_i, \, 3 \delta), \qquad 
U_i= \til{B}(x_i, \, 4 \delta), \qquad i=1, N$$
$$K_1= \text{\rm cl} \big( \til{B}(x_1, \, 3 \delta) \cap K \big)$$
$$K_i= \text{\rm cl} \bigg( \til{B}(x_i, \, 3 \delta) \cap \big( K \setminus
\bigcup_{j=1}^{i-1} K_j \big) \bigg), \qquad i=2,\ldots,N$$
o{\`u} ``cl'' d{\'e}signe l'adh{\'e}rence d'un ensemble. 
Alors $\{ (K_i,\, U_i) \}$ est un d{\'e}coupage de $(K,\, U)$
et nous affirmons qu'il v{\'e}rifie les conditions d{\'e}sir{\'e}es. 

Par construction et hypoth{\`e}se l'application de Gauss est injective
et d'image contenue dans une demi-sph{\`e}re sur chaque $U_i$. De plus, 
$$d_{\Ss^n}(K_i,\ \partial U_i) \ge \delta = \frac1{15}
\bar{d}_{\Ss^n}(K, \, \partial U)  $$
Il ne reste plus qu'{\`a} estimer le nombre de $U_i$ qui peuvent
s'intersecter {\`a} la fois. On utilise le

\begin{lem}\label{boules} 
  Soit $x \in \Ss^n$, $r < \frac\pi2$ et $C \subset B_{\Ss^n}(x, \,
  r)$ un ferm{\'e}. Il existe un recouvrement de $C$ par au plus $2^n$
  boules de rayon $r$ centr{\'e}es en des points de $C$.
\end{lem}

{\it \small Preuve du lemme \ref{boules}:} Soit
$b_n$ le cardinal d'une triangulation de $\Ss^{n-1}=\partial
B_{\Ss^n}(x,\frac\pi2)$ par des simplexes 
de diam{\`e}tre plus petit que $\frac \pi 2$. Pour tout $r < \frac\pi2$
il existe alors une triangulation de $\partial B_{\Ss^n}(x, \,
  r)$ avec $b_n$ simplexes de diam{\`e}tre plus petit que $r$. Les
  pr{\'e}images des simplexes par projection radiale dans $B_{\Ss^n}(x, \,
  r) \setminus \{x\}$ d{\'e}terminent avec $x$ une d{\'e}composition
  $\{S\}$ de
  $B_{\Ss^n}(x, \, r)$ en $b_n$ sous-ensembles de diam{\`e}tre plus petit que
  $r$. Si $C$ intersecte un $S$, alors la boule de rayon $r$ centr{\'e}e
  en un point quelconque de $C\cap S$ contient $S$ et en particulier 
  $C\cap S$. On prend une telle boule pour chaque $S$
  qui intersecte $C$ pour obtenir un recouvrement de $C$ avec au plus
  $b_n$ boules. A titre d'exemple, on peut prendre 
$b_1=2$, $b_2=4$. On prouve alors par r{\'e}curence que $b_n=2^n$
convient.  

\hfill{$\square$}

{\it \small Suite de la preuve du th{\'e}or{\`e}me \ref{theoreme local
    principal} :} 
Deux boules $U_i=\til{B}(x_i, \, 4 \delta)$ et $U_j=\til{B}(x_j, \, 4
\delta)$ peuvent s'intersecter uniquement si $x_j \in \text{\rm cl} \,
\til{B}(x_i, \, 8 \delta)$. Soit $e(n,\, \delta,\, x_i)$ le nombre minimal de
boules de rayon $3\delta$ avec centres dans $\til{B}(x_i, \, 11
\delta)$ n{\'e}cessaires pour recouvrir 
$\text{\rm cl}\, \big( 
\til{B}(x_i, \, 11 \delta) \setminus \til{B}(x_i, \, 4 \delta)
\big)$. 
Alors $U_i$ peut intersecter au plus $2^n \cdot \, e(n,\, \delta,\, x_i)$ parmi
les $U_j$. Dans le cas contraire, avec le Lemme \ref{boules} 
on remplacerait les $B_j$ correspondants 
par {\it au plus} $2^n \cdot \, e(n,\, \delta,\,
    x_i)$ 
   boules de rayon $3\delta$ centr{\'e}es en des points de $K$, sans
toutefois diminuer la partie de $K$ qu'elles recouvrent. 
Mais ceci contredirait la minimalit{\'e} du recouvrement initial. 

La quantit{\'e} $e(n,\, \delta,\, x_i)$ ne d{\'e}pend
pas de $x_i$ car le raisonnement est fait en pratique sur la sph{\`e}re
$\Ss^n$. De plus, on a 
$$e(n,\, \delta,\, x_i) = e(n,\, \delta) \le e_{\RR^n}(n, \, \delta) = 
e(n)$$
o\`u $e_{\RR^n}(n, \, \delta)$ d\'esigne la quantit\'e analogue \`a
$e(n,\, \delta)$ dans l'espace euclidien. Un argument simple
d'homot\'etie montre que $e_{\RR^n}(n, \, \delta)$ ne d\'epend pas de
$\delta$, et elle a \'et\'e not\'ee plus haut par
$e(n)$. L'in\'egalit\'e $e(n,\, \delta) \le e_{\RR^n}(n, \, \delta)$
est assur\'ee par le fait que l'application exponentielle est
contractante sur un espace \`a courbure positive: tout recouvrement
sur $\RR^n$ d\'eterminera par l'exponentielle un recouvrement du
m\^eme cardinal sur la sph\`ere.  

L'existence du d{\'e}coupage 
est donc prouv{\'e}e sous l'hypoth{\`e}se que l'image de 
l'application de Gauss est contenue dans une demi-sph{\`e}re avec les
constantes $C'=15$ et $D'_n=2^n \cdot \, e(n)$. Son existence dans le 
cas g{\'e}n{\'e}ral est prouv{\'e}e 
avec $C=60$ et $D_n=c_n \cdot \, 2^n \cdot \,e(n)$. 
Ceci ach{\`e}ve la d{\'e}monstration.
\hfill{$\square$}

\begin{cor}\label{inegalite presque globale}
 Soit $Z \subset M$ le lieu des points de courbure de Gauss
 nulle. Il existe une constante $C_n $ ind{\'e}pendante de $M$ telle
 que, pour tout compact $K \subset M$ qui n'intersecte pas $Z$ on a 
\begin{equation}\label{maximale} 
\vol(K) \le \bar{C}_n \frac1{\big(\bar{d}_{\: \Ss^n}(K,\,Z)\big)^n} 
 \, \vol \: (Gauss_M)
\end{equation}
\end{cor}

\medskip

\noindent {\bf Remarques.} 
1) Le compact $K$ du corollaire \ref{inegalite presque globale}
peut {\^e}tre disconnexe. 

2)Il est int{\'e}ressant de remarquer le cas $Z=\emptyset$. Ceci
correspond {\`a} une courbure de Gauss partout positive, donc {\`a} une
hypersurface convexe (respectivement {\`a} une courbe localement convexe
dans le cas $n=1$).  
Dans cette situation, le terme de (\ref{maximale}) 
impliquant $\bar{d}$ est absorb{\'e} dans $\bar{C}_n$ 
et nous retrouvons l'in{\'e}galit{\'e} d'Archim{\`e}de 
(section \ref{Applications}) 
avec une constante plus faible et sans ambition de
caract{\'e}riser le cas d'{\'e}galit{\'e}. 

\begin{figure}[h]
 \begin{center} 
  \includegraphics{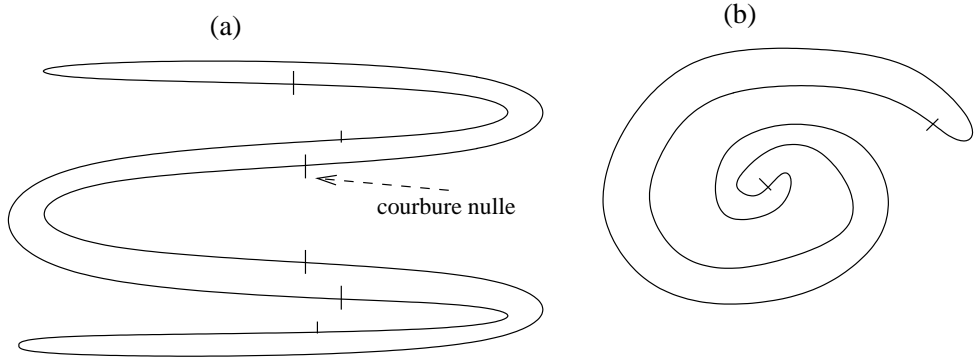} 
 \end{center} 
 \caption{Hypersurfaces avec beaucoup/peu de points de faible
courbure}\label{dessin naif}   
\end{figure}

3) La conclusion de la pr{\'e}sente {\'e}tude est le fait
que le lieu $Z$ des points de courbure nulle joue un r{\^o}le
essentiel dans la comparaison des volumes de $M$ et de son application
de Gauss. L'in{\'e}galit{\'e} (\ref{maximale}) a un caract{\`e}re
asymptotique: plus on voudra attraper dans le compact $K$ 
un volume proche de celui de $M$, plus la distance 
$\bar{d}_{\: \Ss^n}(K,Z)$ devient petite et l'estimation
grossi{\`e}re. L'in{\'e}galit{\'e} est efficace sur des hypersurfaces
enroul{\'e}es sur elles m{\^e}mes (Figure \ref{dessin naif}b), 
mais elle d{\'e}meure tr{\`e}s grossi{\`e}re
pour des hypersurfaces poss{\'e}dant de grandes r{\'e}gions plates ou
ayant un lieu $Z$ trop riche (Figure \ref{dessin naif}a).

\subsection{Techniques globales pour raffiner l'in\'egalit\'e locale}
\label{local to global}

%\setcounter{chapter}{4}
%\setcounter{section}{0}
%\setcounter{equation}{0}
%\section{Local et global s'entrecroisent} \label{local to global} 

%\subsection{Une int{\'e}gration par parties}\label{n n-1}

Reprenons les notations de la section \ref{euler}:
$U\subset M$ est un domaine sur lequel l'application de Gauss est injective
et {\'e}vite au moins un point de $\Ss^n$, $V\subset \Ss^n$, $\Omega
\subset \RR^n$ et $\til{\phi}$ est d{\'e}finie par $\varphi \circ
\til{\phi} = \phi$. 
$$
\xymatrix{
U \ar@/^/[rr]^{\text{Gauss}} \ar@{<-}@/_/[rr]_{\til{\phi}} && V 
\ar@/^/[rr]^{\text{proj.st{\'e}r{\'e}o.}} \ar@{<-}@/_/[rr] && \Omega
}
$$
Ecrivons la formule (\ref{courbure moyenne} ) sur un compact
arbitraire {\`a} bord lisse $K \subset M$: 
\begin{equation} \label{formule utile}
\text{vol}(K)= \frac1n \int_K p\mathcal{H}
dV  + \int_{\partial K} \beta_1
\end{equation} 
avec 
$$\beta_1=\frac{(-1)^{n-1}}{n~!} \ 
 (\varphi,\nu, 
\underbrace{d\varphi,\dots,d\varphi}_{n-1 })$$ 
et 
$$\mathcal{H}=\text{ \rm courbure moyenne sur $M$}$$
Nous allons estimer les termes intervenant dans 
(\ref{formule utile}) 
par les techniques d{\'e}j{\`a} vues dans \ref{Inegalites locales},
sous l'hypoth\`ese $\varphi(M) \subset B^{n+1}(0,R) \subset E^{n+1}$. 
De fa{\c c}on {\'e}vidente on a    
$$\int_{\partial K} \beta_1 \le \frac{R}n \vol(\partial K)$$
Ceci d{\'e}coule de 
$$\vert \beta_1 \vert \le \frac1n \ \text{d} vol _{\partial K}$$
En effet, soit $(f_1,\dots,f_{n-1})$ un rep{\`e}re local orthonorm{\'e} sur
$\partial K$. On a 
$$\vert \beta_1 (f_1,\dots,f_{n-1}) \ \vert \ = \frac1n \
\big| \det \big(\varphi,\nu,\varphi_*(f_1),\dots,\varphi_*(f_{n-1})\big) \big|
\le \frac1n \vert \varphi \vert \le \frac{R}n$$  
puisque $\varphi$ est une isom{\'e}trie et $\varphi(M)$ est contenue dans
la boule de rayon $R$.

Soit maintenant $U \supset K$ un ouvert 
sur lequel l'application de Gauss 
est non-d{\'e}g{\'e}n{\'e}r{\'e}e et injective, avec $\partial U$
lisse.  Nous pouvons alors
exprimer les int{\'e}grales sur $U$ comme des int{\'e}grales sur
$V=Gauss(U)$ et en d{\'e}duire des majorations sur $\int_Kp\mathcal{H}
dV$ en utilisant les
techniques de la section \ref{Inegalites locales}. 
Mettons 
\begin{eqnarray*}
 \alpha_1 & \stackrel{\text{d{\'e}f.}}{=} &  p\mathcal{H}dV  \\
 & = &  \frac{(-1)^n}{(n-1)!} (\varphi, \ d\nu, \
 \underbrace{d\varphi, \dots, \ d\varphi}_{n-1})
\end{eqnarray*}
Evaluons $\til{\phi}^*\alpha_1$ sur un
rep{\`e}re local orthonorm{\'e} $(e_1,\dots,\ e_n)$ 
de $\Ss^n$, qu'on peut choisir comme {\'e}tant constitu{\'e} de directions
de courbure de $M$, c'est-{\`a}-dire $\phi_*e_i=\lambda_ie_i$, avec
$\lambda_i=\frac1{k_i}$ et $k_i$ courbure principale de $M$.
\begin{eqnarray*}
\lefteqn{ \til{\phi}^*\alpha_1 (e_1,\dots,\ e_n)_{\vert_\theta}} \\
& = & \frac{(-1)^n}{(n-1)!} 
 (\varphi, \ d\nu, \ \underbrace{d\varphi, \dots, \ d\varphi}_{n-1}) 
  (\til{\phi}_*e_1,\dots, \ \til{\phi}_* e_n) \\
& = & (-1)^n \sum_{i_1<\dots<i_{n-1}} 
\det\big(f(\theta)\theta+\nabla f(\theta), \ e_1, \dots, \ e_n \big)
\lambda_{i_1} \dots \lambda_{i_{n-1}} \\
& = & (-1)^n f(\theta) \sum_{i_1<\dots<i_{n-1}} 
\lambda_{i_1} \dots \lambda_{i_{n-1}}
\end{eqnarray*}

On a utilis{\'e} $\nu \, \circ \til{\phi} \equiv \text{Id}$, $\varphi \,
\circ \til{\phi} \equiv \phi$. La fonction $f$ est d{\'e}finie par $p \,
\circ \til{\phi} = f$. Dans la somme ci-dessus apparaissent des mineurs de $\phi_*$
(dans la base $(e_1,\dots,\ e_n)$) d'ordre exactement $n-1$. 
Ceci fait qu'ils sont contr{\^o}lables par 
  $\frac{R^{n-1}}{d_{\Ss^n}(K,\ \partial U)^{n-1}}$. Comme $|f| \le R$ 
on obtient
 $$\bigg|  \int_K p \mathcal{H} dV \bigg| \le C_n
\frac{R^n}{d_{\Ss^n}(K,\ \partial U)^{n-1}} \ T(U)$$

L'estimation pr{\'e}c{\'e}dente est valable sous des hypoth{\`e}ses
d'injectivit{\'e} et de bornitude sur l'application de Gauss qui sont
superflues. La technique de globalisation de la section \ref{a bas
  l'injectivite} fournit le r{\'e}sultat suivant: 
\begin{thm}\label{gigigigi}
Soit $\varphi:M^n \longrightarrow E^{n+1}$, $n\ge 2$ 
une immersion isom{\'e}trique d'image contenue dans la boule
$B^{n+1}(0,R)$. 
Soit $U \subset M$ un ouvert sur lequel l'application de Gauss est
non-d{\'e}g{\'e}n{\'e}r{\'e}e. Il existe une constante $\bar{C}_n$
ind{\'e}pendante  
de $U$ et de $M$ telle que, si $K \subset U$ est un compact, on a 
\begin{equation} \label{inegalite de gigigigi}
\vol(K) \le \bar{C}_n \frac{R^n}{\bar{d}_{\Ss^n}(K, \ \partial U)^{n-1}}
T(U) + \frac{R}n \vol(\partial K)
\end{equation}
\end{thm}

La descente de l'exposant de $n$ {\`a} $n-1$ s'appuie sur la formule
int{\'e}grale de Minkowski (\ref{courbure moyenne} ). On pourrait penser
que la formule plus g{\'e}n{\'e}rale (\ref{formule globale generale})
permettrait une descente {\`a} $n-2$. Pourtant ce n'est pas le cas, car
des termes de bord suppl{\'e}mentaires apparaitront et ils ne peuvent
pas {\^e}tre contr{\^o}l{\'e}s de mani{\`e}re raisonnable sauf en dimension
deux. Dans ce dernier cas, (\ref{formule globale generale}) fournit 
$$\vol(K) = \int_K \big( p^2 - \frac12 q^2 \big) \mathcal{K} dS +
\int_{\partial K} \beta _1 + \int_{\partial K} \beta _2$$
avec 
$$\beta_1 = -\frac12 (\varphi, \ \nu, \ d\varphi),$$
$$\beta_2 = -\frac12 p(\varphi, \ \nu, \ d\nu),$$
et  
$$\mathcal{K}= \text{courbure de Gauss sur } M$$
Les in{\'e}galit{\'e}s 
$$ \big| \int_{\partial K} \beta_1 \big| \le \frac{R}2 \Long(\partial K)$$
et 
$$ \big|  \int_K \big( p^2 - \frac12 q^2 \big) \mathcal{K} dS   \big|
\le R^2 T(K)$$
ont d{\'e}j{\`a} {\'e}t{\'e} prouv{\'e}es. Par un calcul similaire {\`a} celui
exprimant l'image inverse de $\beta_1$ on obtient 
$$df \wedge \til{\phi}^*p(\varphi, \ \nu, \ d\nu) = f\vert df \vert ^2
d\vol_{\Ss^n}$$
c'est-{\`a}-dire 
$$\til{\phi}^*\beta_2 = -\frac12 f*df$$
Or $df$ est born{\'e}e par la condition $\varphi(M) \subset B(0,R)$, ce
qui assure que $*df$ est aussi born{\'e}e par $Rd\vol_{Gauss(\partial
  K)}$ {\`a} une constante multiplicative pr{\`e}s. Comme $f$ est elle-m{\^e}me
born{\'e}e par $R$, on d{\'e}duit l'existence d'un $C>0$ tel que  

$$\big|  \int_{\partial K} \beta _2 \big| \le C R^2\cdot \Long(Gauss_{\partial
  K})$$
Il faut remarquer que nous ne supposons pas $Gauss$ comme {\'e}tant
non-d{\'e}g{\'e}n{\'e}r{\'e}e 
sur $K$. Nous en d{\'e}duisons la 

\begin{prop}\label{dim deux} ($n=2$) Soit $\varphi : M^2 \longrightarrow
  \RR^3$ une immersion isom{\'e}trique 
  de surface dont l'image est contenue dans la boule de rayon $R$. 
Soit $U\subset M$ un ouvert {\`a} bord
  lisse. Il existe une constante $C$ telle que   

\begin{equation}\label{isoperimetrique dimension deux}
\Aire(U)\le R^2 \bigg( T(U) + C \, \Long(Gauss_{\partial U}) +
\frac1{2R} \, \Long(\partial U)  \bigg)
\end{equation}
\end{prop}

\noindent {\bf Remarque:} L'in{\'e}galit{\'e} (\ref{isoperimetrique
  dimension deux})  
peut {\^e}tre vue comme une in{\'e}galit{\'e} de
type isop{\'e}rim{\'e}trique: une correction impliquant la fa{\c c}on dont
$M$ est pli{\'e}e est n{\'e}cessaire afin de pouvoir borner l'aire d'un
domaine par la longueur de son bord. 
Par exemple, tout comme on s'y
attend, sur un domaine (presque) plat, l'in{\'e}galit{\'e}
isop{\'e}rim{\'e}trique est (pratiquement) v{\'e}rifi{\'e}e. 

Le caract{\`e}re isop{\'e}rim{\'e}trique est
pr{\'e}sent aussi en dimension sup{\'e}rieure, dans
l'in{\'e}galit{\'e} (\ref{inegalite de gigigigi}), cette fois-ci en
sens inverse: une correction par le volume du bord est n{\'e}cessaire
afin de pouvoir borner de fa{\c c}on optimale le volume d'un domaine
par le volume recouvert sur $\Ss^n$ par son application de Gauss
(l'optimalit{\'e} est entendue au sens de la puissance {\`a} laquelle
apparait $\bar{d}_{\Ss^n}(K,\partial U)$ au d{\'e}nominateur).

%%% Local Variables: 
%%% mode: latex
%%% TeX-master: "Gauss"
%%% TeX-master: "Gauss"
%%% TeX-master: "Gauss"
%%% End: 

%%% Local Variables: 
%%% mode: latex
%%% TeX-master: "Gauss"
%%% End: 

%\input{Gauss4}
%\setcounter{chapter}{5}
%\setcounter{section}{1}
%\setcounter{equation}{0}
\section{Optimalit{\'e} des in{\'e}galit{\'e}s} \label{l'optimalite}

Nous discutons l'optimalit\'e des diff\'erentes in\'egalit\'es que
nous avons obtenues. {\it Pour all\'eger les notations on supposera
d\'esormais que l'image de l'immersion est contenue dans la boule
unit\'e de $E^{n+1}$}. 

\subsection{N{\'e}cessit{\'e} de consid{\'e}rer un {\'e}paisissement
  dans (\ref{locale intrinseque}): }
$$\vol \: (K) \le C_n \frac1{\big(d_{\: \Ss^n}(K,\, \partial U)\big)^n}
 \, \vol \: (Gauss_U) $$
Le membre de droite est optimal au sens o{\`u} on ne peut remplacer
$Gauss_U$ par $Gauss_K$. Sous les hypoth{\`e}ses du corollaire \ref{coro}
une in{\'e}galit{\'e} du type 
\begin{equation}\label{U remplace par K}
\vol \: (K) \le C_n \frac1{d^n}
 \, \vol \: (Gauss_K) 
\end{equation} 
ne peut {\^e}tre vraie pour toute hypersurface $M$ et toute paire $(U,
\, K)$ telle que $d_{\Ss^n}(K, \, \partial U)
\ge d $ fix{\'e}, comme le montre l'exemple suivant. 

Pour plus d'aisance dans l'{\'e}criture, on fait la construction en
dimension $1$. On consid{\`e}re la famille d'ellipses 
$$E_\epsilon=\big\{ (x,y)~: x^2+\frac{y^2}{\epsilon^2}=1 \big\}
\subset B(0,1)$$ 
On fixe $\frac\pi2 > d > 0$. On choisit 
$U_\epsilon =\{ (x,y) \in E_\epsilon~: \, y>0 \}$ la demi-ellipse
sup{\'e}rieure et $K_\epsilon = \{ (x,y)\in U_\epsilon~: \vert x \vert \le
r\} $, avec $0 < r < 1$ fix{\'e}. 
Pour $\epsilon$ suffisamment petit on a effectivement
$d_{\Ss^1}(K_\epsilon, \, U_\epsilon) \ge d$, mais
$\Long(K_\epsilon)\tto{\epsilon\rightarrow 0} 2r$, tandis que
$\Long(Gauss_{K_\epsilon})\tto{\epsilon\rightarrow 0} 0$. Ceci montre
qu'une in{\'e}galit{\'e} du type (\ref{U remplace par K}) ne peut {\^e}tre
v{\'e}rifi{\'e}e. 
La m{\^e}me construction fonctionne en dimension sup{\'e}rieure en
consid{\'e}rant la famille d'ellipso{\"\i}des 
$$E_\epsilon=\big\{ (x_1,\dots,x_{n+1}) \in \RR^{n+1}~: x_1^2+\dots+x_n^2+\frac{x_{n+1}^2}{\epsilon^2}=1 \big\}
\subset B(0,1)$$ 
ainsi que 
$$U_\epsilon= \big\{ (x_1,\dots,x_{n+1}) \in E_\epsilon~: \, x_{n+1}>0
\big\}, \qquad
K_\epsilon = \big\{ (x_1,\dots,x_{n+1})\in U_\epsilon~: \sum_{i=1}^n x_i
^2\le r^2 \big\}$$

\subsection{Optimalit{\'e} de l'in{\'e}galit{\'e}
(\ref{isoperimetrique dimension deux})} 
Les trois termes
$\Aire(Gauss_U)$, $\Long(Gauss_{\partial U})$ et $\Long(\partial U)$
constituant la partie droite de l'in{\'e}galit{\'e} 
sont ind{\'e}pendants, dans le sens qu'aucun d'entre eux ne peut
{\^e}tre born{\'e} par les deux autres. Dans la figure \ref{exemples 
figuratifs pour l'inegalite isoperimetrique en dimension 2} nous
donnons trois exemples de bouts de
surfaces {\`a} bord dans lesquels deux des trois termes sont
n{\'e}gligeables par rapport au troisi{\`e}me. 
Les termes n{\'e}gligeables sont, respectivement: $\Aire(Gauss_U)$ et
$\Long(Gauss_{\partial U})$ pour (1); $\Long(\partial U)$ et
$\Long(Gauss_{\partial U})$ pour (2); $\Aire(Gauss_U)$ et
$\Long(\partial U)$ pour (3). 

\begin{figure}[h]  
%\vspace{0.1cm} 
\begin{center} 
\includegraphics*{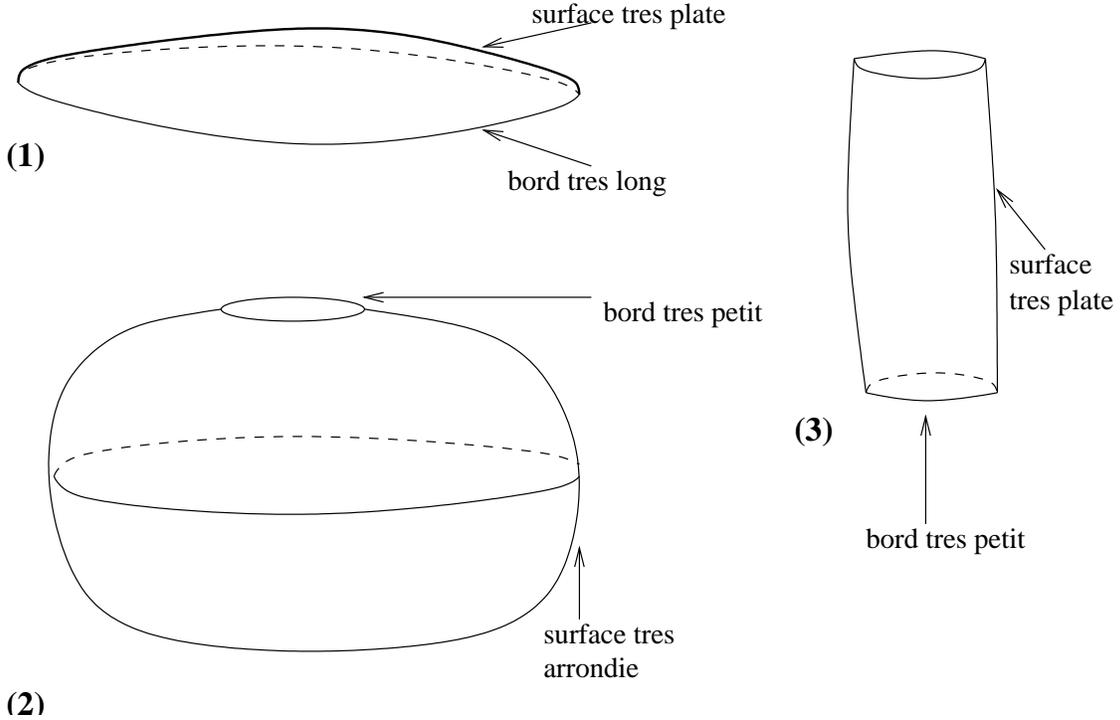} 
\end{center} 
%\vspace{6cm} 
   \caption{Les quantit{\'e}s $\Aire(Gauss_U)$, $\Long(Gauss_{\partial
       U})$ et $\Long(\partial U)$ sont ind{\'e}pendantes}
       \label{exemples figuratifs pour l'inegalite isoperimetrique en
         dimension 2}   
\end{figure}

Les exemples $(1)$ et $(2)$ dans la figure \ref{exemples 
figuratifs pour l'inegalite isoperimetrique en dimension 2} prouvent
aussi que les quantit{\'e}s $\text{Aire}(Gauss_U)$ et
$\text{Long}(\partial U)$ sont essentielles dans le fonctionnement de
l'in{\'e}galit{\'e} (\ref{isoperimetrique dimension deux}). 
On n'a pas encore d'exemple pour montrer que la
quantit{\'e} $\text{Long}(Gauss_{\partial U})$ est aussi essentielle. Il
est fort possible qu'une in{\'e}galit{\'e} du type
$$\text{Aire}(U) \le C\text{Aire}(Gauss_U) + C'\text{Long}(\partial U)$$
soit v{\'e}rifi{\'e}e.

\subsection{Meilleure estimation pour le volume d'un compact $K
  \subset M$ fix{\'e}} \label{belle motivation}
Pla{\c c}ons-nous dans le contexte du corollaire
  \ref{coro}: $K$ est un compact admettant un voisinage $U$ sur lequel
  $Gauss$ est non-d{\'e}g{\'e}n{\'e}r{\'e}e, injective et d'image contenue
  dans une demi-sph{\`e}re. L'in{\'e}galit{\'e} (\ref{locale intrinseque})
  est alors valable:
  $$\vol \: (K) \le C_n \frac1{\big(\text{d}_{\: \Ss^n}(K,\, \partial
    U)\big)^n}  \, \vol \: (Gauss_U) $$
  Le but initial de cette {\'e}tude {\'e}tait de majorer le volume de $K$
  par un certain volume recouvert sur $\Ss^n$ par $Gauss$. Il est
  {\'e}vident que, pour une distance $d = d_{\Ss^n}(K, \ \partial U)$
  {\it fix{\'e}e}, la meilleure estimation par (\ref{locale intrinseque}) est
  obtenue en choissant $U$ tel que $Gauss(U)=Gauss(K)_d$, o{\`u} 
  $Gauss(K)_d $ est l'{\'e}paisissement par $d$ de $Gauss(K)$:
  $$Gauss(K)_d = \big\{ x\in \Ss^n \ : \ \text{dist}_{\Ss^n}(x, \ Gauss(K)) <
   d  \big\}$$
  On {\'e}crira $K_d$ pour l'{\'e}paisissement par $d$ de $K$ par rapport
  {\`a} la distance sph{\'e}rique sur $U$: 
  $$K_d= \big\{ x\in U \ : \ d_{\Ss^n}(x, \ K) <
   d  \big\}$$
  Alors $Gauss(K_d)= Gauss(K)_d$. 

  Ainsi $K_d$ est le plus petit ouvert $\mathcal{U}$ satisfaisant
  $d_{\Ss^n}(K, \ \partial \mathcal{U}) = d $ et on a, sous l'hypoth{\`e}se
  d'existence d'un $U$ comme ci-dessus,
  \begin{equation}\label{epaisissements}
   \vol(K) \le C_n \frac 1 {d^n} \vol (Gauss(K)_d)
  \end{equation}
  Cette {\'e}criture de l'in{\'e}galit{\'e} ne fait intervenir de voisinage
  du compact 
  qu'implicitement. Le $d$ peut varier entre $0$ et $d_{\max}$, o{\`u} 
  \begin{equation*} \left. \begin{array}{l}
  d_{\max} = \max \big\{ d \ : \ \exists\ U \text{ ouvert }
  \supset K_d \text{ t.q. }
  Gauss_{\vert_U} \text{ non-d{\'e}g{\'e}n{\'e}r{\'e}e,}    \\ 
  \qquad \qquad \qquad  \text{ injective et d'image
    contenue dans une demi-sph{\`e}re } \big\} \end{array} \right.
  \end{equation*}

  Nous nous int{\'e}ressons alors {\`a} comprendre quel est le minimum du
  membre de droite dans (\ref{epaisissements}). Le probl{\`e}me plus
  g{\'e}n{\'e}ral {\`a} traiter est de trouver 
  $$\min_{0< d < d_{\max}} \frac{\vol(A_d)}{d^n}$$
  pour $A\subset \Ss^n$ un compact quelconque inclus dans une
  demi-sph{\`e}re $\mathcal{D}$, $A_d$ un {\'e}paisissement par $d$ de $A$
  et $d_{\max} \le \text{dist}_{\Ss^n}(A,\
  \partial \mathcal{D})$. La r{\'e}ponse compl{\`e}te {\`a} cette question
  est fournie par la proposition suivante, d{\'e}montr{\'e}e dans un cadre
  plus g{\'e}n{\'e}ral dans l'appendice. 

  \begin{prop} \label{la plus belle}
   Pour tout compact $K \subset \Ss^n$ inclus dans une demi-sph{\`e}re
   la fonction 
   $$f(d)  = \frac{\vol(Gauss(K)_d)}{d^n}$$
   est d{\'e}croissante sur l'intervalle $]0,d_{\max}]$. 
  \end{prop} 

  En particulier, la meilleure estimation dans l'in{\'e}galit{\'e}
  (\ref{epaisissements}) est obtenue pour $d=d_{\max}$.

\subsection{Optimalit{\'e} de l'in{\'e}galit{\'e} (\ref{gener}) } 
\begin{equation} \label{label local}
\vol(K) \le \bar{C}_n \frac 1 {\bar{d}_{\Ss^n}(K, \ \partial U)^{n-1} }
\vol(Gauss_U) + \frac 1n \vol( \ \partial K) 
\end{equation}

Nous rapellons que $Gauss$ est suppos{\'e}e non-d{\'e}g{\'e}n{\'e}r{\'e}e sur
l'ouvert $U$ et $K \subset U$ est un compact. De plus, l'in{\'e}galit{\'e}
fonctionne pour une {\it immersion} isom{\'e}trique de $U$ dans
$\RR^{n+1}$.  
Nous allons montrer que les deux termes de droite dans (\ref{label local}) 
sont essentiels dans le fonctionnement de l'in{\'e}galit{\'e} et que
l'exposant $n-1$ est optimal.  

a) Le terme $\vol(\ \partial K)$ est essentiel. Consid{\'e}rons une
suite $U_k$ de bouts de surface d{\'e}finis par
\begin{equation} \label{hemispheres}
U_k = B^{n+1}(0,1) \cap S^n\big((0,\dots,\ 0, \ -k), k+ \frac 1k
\big)
\end{equation} 
Les $U_k$ sont des h{\'e}mi-sph{\`e}res de courbure de plus en plus
faible, qui tendent avec $k$ vers le disque ouvert
$$D^n(0,1) = \{ x= (x_1,\dots,\ x_{n+1}) \in \RR^{n+1} \ : \  \vert x
\vert < 1 \text{ et } x_{n+1}=0 \} $$
L'image par $Gauss$ de $U_k$ est une boule g{\'e}od{\'e}sique sur $\Ss^n$
centr{\'e}e au p{\^o}le nord,
de rayon $\rho_k \rightarrow 0$. Soit $K_k$ la pr{\'e}image de la boule
g{\'e}od{\'e}sique ferm{\'e}e de m{\^e}me centre et rayon $\rho_k/2$. Alors
$K_k \rightarrow \bar{D}^n(0,\ \frac 12)$ (Figure \ref{l'inegalite est
  optimale}, (1)). 
Par construction $d_{\Ss^n}(K_k,\ \partial U_k) = \rho_k/2$, tandis
que $\vol(Gauss_{U_k}) \sim \rho_k ^n$. Ainsi 
$$ \frac
{\vol(Gauss_{U_k})} {\big(d_{\Ss^n}(K_k,\ \partial U_k) \big) ^{n-1}}
\sim \rho_k \longrightarrow 0$$
devient n{\'e}gligeable devant $\vol(K_k)$. 
Ceci prouve que le terme $\vol(\ \partial K_k )$ est essentiel pour le
fonctionnement de l'in{\'e}galit{\'e} 
(\ref{gener}). Le m{\^e}me exemple prouve que
{\bf l'in{\'e}galit{\'e} (\ref{gener}) est optimale} au sens o{\`u} l'exposant $n$
ne peut {\^e}tre diminu{\'e} sans ajout de termes suppl{\'e}mentaires. 

\medskip 

b) Le terme $\frac 1 {\bar{d}_{\Ss^n}(K, \ \partial U)^{n-1} }
\vol(Gauss_U)$ est essentiel. Une in{\'e}galit{\'e} du type 
$\vol(K) \le C \vol(\ \partial K)$
ne peut {\^e}tre vraie telle quelle, comme on peut le voir facilement
en prenant pour $K$ une sph{\`e}re dont on a enlev{\'e} une boule
g{\'e}od{\'e}sique arbitrairement petite. 

\medskip 

c) L'exposant $n-1$ dans (\ref{gener}) est
optimal pour le coefficient $\frac 1n$ devant $\vol(\ \partial K)$. 
A l'int{\'e}rieur du disque plan 
$$D^2(0,1) = \{ x= (x_1,\dots,\ x_{n+1}) \in \RR^{n+1} \ : \  \vert x
\vert < 1 \text{ et } x_3 = \dots = x_{n+1}=0 \} $$
consid{\'e}rons une courbe ferm{\'e}e imerg{\'e}e,
appel{\'e}e $\gamma$. 

La construction que nous pr{\'e}sentons g{\'e}n{\'e}ralise celle des
surfaces de rotation dans $\RR^3$. La courbe $\gamma$ admet en chaque
point un plan normal g{\'e}om{\'e}trique de dimension $n$, et tous ces plans
normaux sont canoniquement identifi{\'e}s par des translations et
rotations dans $\RR^2 \equiv \{ x_3 = \ldots  = x_{n+1} = 0\}$. Fixons une
hypersurface $H$ 
de dimension $n-1$ passant par l'origine dans un de ces
plan normaux. Son image via les identifications canoniques d{\'e}crit le
long de $\gamma$ une hypersurface $\mathcal{H}$ de dimension $n$ 
dans $\RR^{n+1}$
qu'on appelle {\it hypersurface de rotation} d{\'e}finie par $H$ le
long de $\gamma$. 

Nous donnons des exemples de telles hypersurfaces de rotation
$\mathcal{H}$ pour lesquelles la diff{\'e}rence $\vol(K)- \frac 1n
\vol(\ \partial K)$ est arbitrairement grande et domin{\'e}e par
$C_n \frac 1 {\bar{d}_{\Ss^n}(K, \ \partial U)^{n-1}} \vol(Gauss_U)$,
mais pour lequelles $C \frac 1 {\bar{d}_{\Ss^n}(K, \ \partial
  U)^{n-2}} \vol(Gauss_U)$ est arbitrairement petit quelle que soit la
constante $C$. 

Faisons d'abord une remarque concernant la lissit{\'e} des hypersurfaces
$\mathcal{H}$ construites comme ci-dessus. Comme le montre l'exemple
des sph{\`e}res (pour lesquelles on prend pour $\gamma $ un cercle et
pour $H$ une demi-sph{\`e}re), $\mathcal{H}$ peut avoir des
singularit{\'e}s, issues de l'intersection de "m{\'e}ridiens" $H$ voisins.
L'absence de singularit{\'e}s sur un m{\'e}ridien
$H_p$ au point $p\in \gamma$ d{\'e}pend de la courbure de $\gamma$ en $p$
et de la courbure de $H$. En termes pr\'ecis, il faut que la projection
de $H$ sur la normale \`a $\gamma$ au point $p$ s'\'eloigne de $p$
d'une distance au plus \'egale au rayon focal en $p$. Comme cela on
est sur que les $H$ ``voisines'' ne s'intersectent pas. 
En particulier, si la courbure de $\gamma$ est born{\'e}e
sup{\'e}rieurement, on peut prendre pour $H$ une h{\'e}misph{\`e}re
centr{\'e}e sur 
$\gamma$, de courbure suffisamment faible. 
Ceci est le point de d{\'e}part de notre construction.

\medskip

Prenons donc pour $\gamma$ une courbe ferm{\'e}e imerg{\'e}e
localement convexe. Pour tout 
$\delta < 1$ il existe une telle courbe de longueur
{\it arbitrairement grande} contenue dans $D^2(0,\delta)$ et de
courbure born{\'e}e sup{\'e}rieurement par $2/\delta$. En particulier
le rayon focal est minor\'e par $\delta/2$. 
Soit $\mathcal{H}_k$ une suite d'hypersurfaces de rotation
construites sur $\gamma$ {\`a} partir d'h{\'e}misph{\`e}res de rayon
g{\'e}od{\'e}sique {\'e}gal {\`a} $r<1$ et courbure tendant vers $0$ avec $k$. Le $r$
est choisi pour que $\mathcal{H}_k \subset B^{n+1}(0,1)$. En prenant
$\delta $ suffisamment petit, le $r$ peut {\^e}tre choisi {\it
  arbitrairement proche de $1$}. 

\begin{figure}[h]  
\begin{center} 
\includegraphics{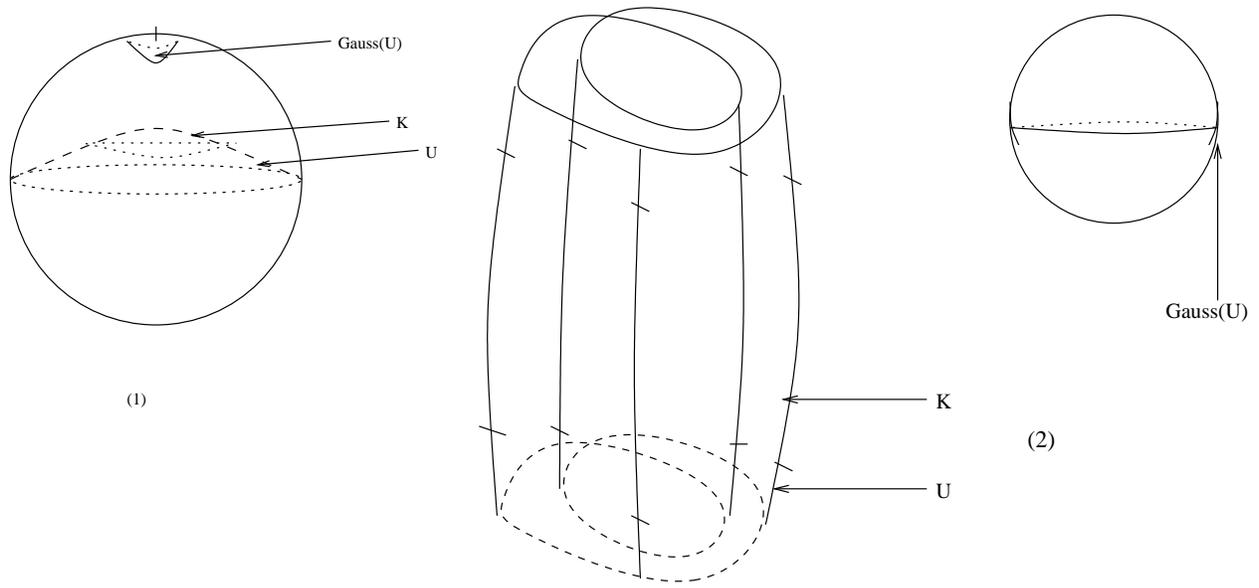} 
\end{center} 
   \caption{Optimalit{\'e} de l'in{\'e}galit{\'e} (\ref{gener})
       \label{l'inegalite est optimale}}   
\end{figure}

Prenons $U_k=\mathcal{H}_k$. L'image par Gauss de $U_k$ est un tube de
rayon $\rho_k \longrightarrow 0$ autour du grand cercle
$\bar{D}^2(0,1) \cap \Ss^n$, et on d{\'e}finit $K_k$ comme la pr{\'e}image
par Gauss du tube de rayon $\rho_k/\lambda$ autour du m{\^e}me grand
cercle, avec $\lambda > 1$ (Figure \ref{l'inegalite est optimale},
(2)). Alors 
$$\vol(K_k) \longrightarrow \vol \ B^{n-1}(\frac r \lambda) \cdot
\Long(\gamma) = \frac r \lambda \frac 1 {n-1} \vol \ \Ss^{n-2}( \frac
r \lambda ) \cdot \Long(\gamma)$$ 
$$\vol(\ \partial K_k) \longrightarrow \vol \ \Ss^{n-2}( \frac
r \lambda ) \cdot \Long(\gamma)$$ 
Si $\frac r \lambda > \frac {n-1} n$ (ce qui est r{\'e}alisable pour
$r$, $\lambda$ proches de $1$), alors $\vol(K_k) - \frac 1n \vol (\
\partial K_k ) $ est positif et peut {\^e}tre rendu arbitrairement grand
en augmentant $\Long(\gamma)$. Ceci d\'emontre d\'ej\`a qu'une
correction de $\vol(K_k)$ par un terme suppl{\'e}mentaire est
n{\'e}cessaire. 
D'un autre c{\^o}t{\'e}, 
$$ \vol (Gauss_{U_k}) \sim \rho_k ^{n-1} \cdot \Long(\gamma)$$
$$ \bar{d}_{\Ss^n}(K_k,\ \partial U_k) = \rho_k \big(  1 - \frac 1
\lambda \big)$$
Ceci montre qu'un exposant $\eta$ strictement plus petit que $n-1$ rendrait
la quantit{\'e} $1 / {\big( \bar{d}_{\Ss^n}(K_k, \ \partial U_k)
  \big)^\eta}  \cdot 
\vol(Gauss_{U_k})$ {\'e}quivalente {\`a} $\rho_k^{n-1-\eta}
\longrightarrow 0$, ne pouvant donc pas dominer $\vol(K_k) - \frac 1n \vol (\
\partial K_k ) $. L'exposant $n-1$ est donc optimal pour le
coefficient $\frac 1 n $ dans (\ref{gener}).

\hfill{$\square$}

{\it Remarque.} Nous ne savons pas prouver que l'exposant $n-1$ est
optimal pour tout coefficient devant $\vol(\partial K)$. Par exemple,
notre exemple ne fonctionne pas si le coefficient $\frac1n$ est remplac\'e par 
$\frac{1}{n-1}$.

%%% Local Variables: 
%%% mode: latex
%%% TeX-master: "Gauss"
%%% End: 

%\setcounter{chapter}{7}
%\setcounter{section}{1}
%\setcounter{equation}{0}
\section{Appendice I. Une in{\'e}galit{\'e} isop{\'e}rim{\'e}trique inverse}

  Les r\'esultats que nous pr\'esentons dans cet appendice ont \'et\'e
motiv\'es par la recherche d'une preuve pour la proposition \ref{la
plus belle}, \`a son tour motiv\'ee par l'\'etude de la meilleure
fa\c{c}on d'estimer le volume d'un compact \ref{belle motivation}. 
  Dans la suite on d\'esignera par $\mc{R}^n$ l'espace
  euclidien $\RR^n$, la sph{\`e}re $\Ss^n$ ou l'espace hyperbolique
  $H^n$. Si $K$ est un compact de $\mc{R}^n$ et $d>0$ un r{\'e}el
  positif on note 
  $$K_d = \{ \ x \in \mc{R}^n \ : \ d(x, \ K) < d \ \}$$
  On note par $B(d)$ une boule ouverte de rayon $d$ dans $\mc{R}^n$. 
  On pose
  $$d_{\tx{max}}(K) \ = \ \pi - \tx{Circumradius}(K), \qquad \qquad \qquad
  \mc{R}^n = \Ss^n$$ 
  $$  d_{\tx{max}}(K) = + \infty, \qquad \qquad \qquad \qquad \qquad \qquad
   \mc{R}^n = \RR^n, \ H^n$$ 
  Nous allons montrer la 
  \begin{prop}\label{isoperimetrique inverse}
    Pour tout compact $K \subset \mc{R}^n$ et pour tout r{\'e}el $0 < d
    < d_{\tx{max}}(K)$ on a l'in{\'e}galit{\'e} 
    \begin{equation}\label{inegalite isoperimetrique inverse} 
      \frac{\vol_n \ K_d}{\vol_{n-1}\ \partial K_d} \ge \frac{\vol_n \
      B(d)}{\vol_{n-1} \ \partial B(d)}
    \end{equation}
    avec {\'e}galit{\'e} si et seulement si $K$ est un ensemble fini
    de points situ{\'e}s {\`a} distance au moins $d$ l'un de
    l'autre (ou, en d'autres mots, si $K_d$ est une union disjointe de
    boules de rayon $d$). 
  \end{prop}

  Avant d'en donner la preuve, 
  indiquons un corollaire et une formulation {\'e}quivalente de la proposition
  \ref{isoperimetrique inverse}.  
  \begin{cor}\label{corollaire interessant}
   Sous les hypoth{\`e}ses pr{\'e}c{\'e}dentes et pour $\mc{R}^n= \RR^n, \
   \Ss^n$ on a 
   $$\frac{\vol_n \ K_d}{\vol_{n-1}\ \partial K_d} \ge \frac dn $$
   En particulier, la fonction 
   \begin{equation}\label{quotient} 
     f(d) = \frac{\vol_n \ K_d}{d^n}
   \end{equation} 
   est d{\'e}croissante pour $0 < d < d_{{\tx{max}}}(K)$. 
  \end{cor}

  \demo 
 Soit $\alpha_{n-1}$ le volume de la sph{\`e}re unit{\'e} de
  $\RR^n$. Il est facile de voir que~: 
  \begin{itemize} 
   \item pour $\mc{R}^n = \RR^n$, on a $\vol_{n-1} \ \partial B(d) =
     d^{n-1} \alpha_{n-1}$, $\vol_n \ B(d) = \int_0^d t^{n-1}
     \alpha_{n-1} dt = \frac {d^n}{n} \alpha_{n-1}$. Ceci
     entra{\^\i}ne 
     $$\frac{\vol_n \ B(d) }{\vol_{n-1} \ \partial B(d) } =
     \frac{d}{n}$$
   \item pour $\mc{R}^n = \Ss^n$, on a  $\vol_{n-1} \ \partial B(d) =
  (\sin d)^{n-1} \alpha_{n-1}$, $\vol_n \ B(d) = \int_0^d (\sin
  t)^{n-1} \alpha_{n-1} dt$. Alors on montre que 
  $$\frac{\vol_n \ B(d) }{\vol_{n-1} \ \partial B(d) } = \frac{\int_0^d (\sin
  t)^{n-1} dt}{(\sin t)^{n-1}} \ge \frac dn$$
   Pour cela il suffit de voir que la fonction $g(d) = n
   \int_0^d (\sin t)^{n-1} dt - d (\sin d)^{n-1} $ est croissante et
   nulle en $d=0$, ce qui assure sa positivit{\'e}. 
  \end{itemize} 
  \hfill{$\square$}

  \medskip

  {\bf Remarque.} C'est l'{\'e}tude de la fonction $f(d)$ dans
  (\ref{quotient}) qui a engendr{\'e} l'{\'e}tude de (\ref{inegalite
  isoperimetrique inverse}). Dans le cas d'un convexe $K \subset
  \RR^n$ une formule due {\`a} Steiner (\cite{BZ}, \S 4.19) assure que
  $\vol_n \ K_d$ est un polyn{\^o}me de degr{\'e} $n$ en $d>0$ ayant les 
  coefficients positifs, d'o{\`u} la d{\'e}croissance de $f$. Dans le cas
  non-convexe on peut encore {\'e}crire $\vol_n \ K_d$ comme polyn{\^o}me de
  degr{\'e} $n$ en $d$ \cite{Fe}, mais uniquement pour des valeurs de
  $d$ plus petites que le rayon focal de $\partial K$ lorsque celui-ci
  est lisse, ou plus petites qu'un analogue du rayon focal pour des
  compacts {\`a} bord "raisonnablement" singulier \cite{Fe}. 

  C'est la d{\'e}croissance de $f$ pour de tr{\`e}s petites et de tr{\`e}s
  grandes valeurs de $d$ (dans ce dernier cas, $K_d$ est proche d'une
  boule de rayon $d$) qui nous a men{\'e}s vers l'{\'e}tude de sa
  d{\'e}croissance globale.

  Il convient d'observer que la variation de $\vol_n \ K_d$ peut 
  {\^e}tre tr{\`e}s compliqu{\'e}e dans le domaine interm{\'e}diaire, 
  notamment {\`a} cause de la pr{\'e}sence de
  trous dans $K$ qui peuvent {\^e}tre absorb{\'e}s dans $K_d$ pour $d$
  suffisamment grand. 
  
  \bigskip

  La classe d'ouverts de $\mc{R}^n$ pour
  lesquels l'in{\'e}galit{\'e} \ref{inegalite isoperimetrique inverse} est
  valable est 
  $$\mc{C}_d = \big\{ D \text{ \rm ouvert born{\'e}} \subset \mc{R}^n \ : \
  \forall \ x \in D, \ \exists \text{ \rm boule } x\in B(d) \subset D \
  \big\}$$ 
  De fa{\c c}on {\'e}vidente, si $K$ est un compact alors 
  \begin{equation}\label{inclusion}
   K_d \in \mc{C}_d
  \end{equation}
  La proposition \ref{isoperimetrique inverse} d{\'e}coulera
  alors de la 
  \begin{prop}\label{isoperimetrique inverse equivalente}
   Pour tout ouvert $D \in \mc{C}_d$ on a 
    \begin{equation}\label{inegalite a prouver}
      \frac{\vol_n \ D}{\vol_{n-1} \ \partial D} \ge \frac {\vol_n \
      B(d)}{\vol_{n-1} \ \partial B(d)}
    \end{equation} 
   On a {\'e}galit{\'e} si et seulement si $D$ est une union disjointe de
   boules de rayon $d$. 
  \end{prop}
  C'est cette derni{\`e}re in{\'e}galit{\'e} qui peut {\^e}tre regard{\'e}e
  comme une in{\'e}galit{\'e} isop{\'e}rim{\'e}trique ``inverse''.  Remarquons
  d'ailleurs que $\mc{C}_d$ co{\"\i}ncide avec 
  $$\mc{C}'_d = \big\{ D \text{ \rm ouvert} \subset \mc{R}^n \ : \
  \exists \ K \text{ \rm compact, } D = K_d \big\}$$
  En effet, l'inclusion $\mc{C}'_d \subseteq \mc{C}_d$ d{\'e}coule de
  (\ref{inclusion}). D'un autre c{\^o}t{\'e} tout ouvert $D \in \mc{C}_d$
  peut s'{\'e}crire
  $$D = (D_{-d})_d$$
  avec 
  $$D_{-d} = \{ x \in D \ : \ d(x, \ \partial D ) \ge d  \}$$
  L'inclusion $(D_{-d})_d \subseteq D$ est valable pour un ouvert
  arbitraire, n'appartenant pas n{\'e}cessairement {\`a} $\mc{C}_d$. 
  Pour voir l'inclusion inverse, prenons un point $x \in D$
  quelconque. Puisque $D \in \mc{C}_d$ il existe une boule $x \in
  B(y,d) \subset D$, ce qui assure que $y\in D_{-d}$. Ceci
  entra{\^\i}ne 
  $B(y,d) \in (D_{-d})_d$ et donc $x \in (D_{-d})_d$. 

  On conclut que les propositions
  \ref{isoperimetrique inverse} et 
  \ref{isoperimetrique inverse equivalente} sont {\'e}quivalentes. 

  \medskip 

\begin{figure}[h]  
\begin{center} 
\includegraphics{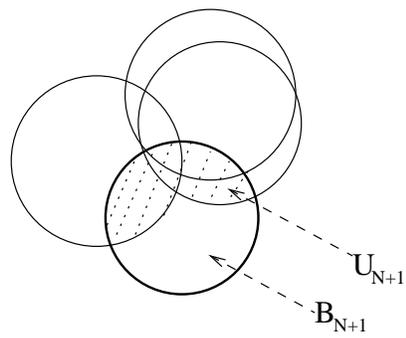} 
\end{center} 
   \caption{Preuve de la proposition \ref{isoperimetrique inverse equivalente}
       \label{Boules}}   
\end{figure}

  {\small \it D{\'e}monstration de la proposition \ref{isoperimetrique
  inverse equivalente}.} \begin{footnote}
{Une premi\`ere preuve, en dimension $2$ et pour la m\'etrique euclidienne, 
m'a \'et\'e sugg\'er\'ee par Nicolae Mihalache. Je lui en suis
vivement reconnaissant.  
}
\end{footnote}
Un ouvert $D \in \mc{C}_d$ peut s'{\'e}crire
  comme
  $$D = \bigcup_{i=1}^\infty B_i(d)$$
  avec 
  $$ \lim _N \vol_n \ \bigcup _{i=1} ^N B_i(d)  = \vol_n \ D$$
  $$ \lim _N \vol_{n-1} \ \partial \bigcup _{i=1} ^N B_i(d)  = \vol_{n-1}
  \ \partial D$$
  Il suffit donc de prouver (\ref{inegalite a prouver}) lorsque $D$
  est une union finie de boules de rayon $d$. On raisonne par
  r{\'e}currence sur $N$. Lorsque $N=1$ on a m{\^e}me l'{\'e}galit{\'e}. Pour
  all{\'e}ger les notations dans le passage de $N$ {\`a} $N+1$ on note
  $$V_N = \vol_n \ \bigcup _{i=1} ^N B_i(d) \qquad \qquad \qquad  
  \partial _N = \vol_{n-1} \ \partial \bigcup _{i=1} ^N B_i(d)  $$
  Soit 
  $$U_{N+1} = B_{N+1}(d) \cap \Big( \bigcup _{i=1} ^N B_i(d)  \Big)$$
  Alors (voir aussi la figure \ref{Boules}) 
  $$V_{N+1} + \vol _n \ U_{N+1} = V_N + \vol_n \ B_{N+1}(d) $$
  $$\partial _{N+1} + \vol_{n-1} \ \partial U_{N+1} = \partial _N +
  \vol_{n-1} \ \partial B_{N+1} (d) $$
  En appliquant l'hypoth{\`e}se de r{\'e}currence on d{\'e}duit
  $$\frac{V_{N+1} + \vol _n \ U_{N+1} }{\partial _{N+1} + \vol_{n-1} \
    \partial U_{N+1} } \ge \frac{\vol_n \ B(d) }{\vol _{n-1} \
    \partial B(d)}$$
  Pour conclure il suffit alors de montrer que 
  \begin{equation}\label{inegalite intermediaire}
    \frac{\vol _n \ U_{N+1} }{\vol_{n-1} \
    \partial U_{N+1}} \le \frac{\vol_n \ B(d) }{\vol _{n-1} \
    \partial B(d)}
  \end{equation} 
  On utilise maintenant l'in{\'e}galit{\'e} isop{\'e}rim{\'e}trique dans
  $\mc{R}^n$ (\cite{BZ}, \S 10.2.1): 
  si $B(h)$ d{\'e}signe la boule de volume {\'e}gal {\`a}
  $\vol_n \ U_{N+1}$ on a l'in{\'e}galit{\'e} $\vol_{n-1} \ \partial
  U_{N+1} \ge \vol_{n-1} \ \partial B(h)$. On en d{\'e}duit 
  $$\frac{\vol _n \ U_{N+1} }{\vol_{n-1} \
    \partial U_{N+1}} \le \frac{\vol_n \ B(h) }{\vol _{n-1} \
    \partial B(h)}$$
  D'un autre c{\^o}t{\'e}, puisque $U_{N+1} \subseteq B_{N+1}(d)$ on a
  forc{\'e}ment $h \le d$ et l'in{\'e}galit{\'e} (\ref{inegalite
    intermediaire}) s'ensuit une fois qu'on aura remarqu{\'e} que la
  fonction 
  $G(d) = \frac{\vol_n \ B(d)}{\vol_{n-1} \ \partial B(d) }$
  est strictement croissante en $d>0$.  
  \begin{itemize} 
   \item si $\mc{R}^n = \RR^n$ alors $G(d) = \frac dn$~;
   \item si $\mc{R}^n = \Ss^n$ alors $G(d) = \frac {\int_0^d (\sin
     t)^{n-1} dt }{(\sin t)^{n-1}} $. Sa croissance {\'e}quivaut {\`a} la
   positivit{\'e} de $G_1(d)  = (\sin d)^n - (n-1) \int_0^d(\sin
   t)^{n-1} dt \cos d$, or $G_1(0) = 0$ et $G_1'\ge 0$~;
   \item si $\mc{R}^n = H^n$ alors $G(d) = \frac {\int_0^d (\tx{sinh }
     t)^{n-1} dt }{(\tx{sinh } t)^{n-1}} $. Sa croissance {\'e}quivaut
   {\`a} la positivit{\'e} de $G_1(d) = (\sinh d)^n - (n-1)\int_0^d(\sinh
   t)^{n-1} dt \cosh d$, or $G_1(0)=G_1'(0) =0$ et $G_1'' \ge 0$.
  \end{itemize} 

 \hfill{$\square$}

%%% Local Variables: 
%%% mode: latex
%%% TeX-master: "Gauss"
%%% End: 

\end{document}